\documentclass[11pt]{olplainarticle}

\usepackage{amsmath}
\usepackage{graphicx}
\usepackage{subcaption}
\usepackage{algorithm}
\usepackage{algpseudocode}
\usepackage{mypackage}
\usepackage{mathtools}
\mathtoolsset{showonlyrefs}

\title{Total variation regularization with reduced basis in electrical impedance tomography}

\author[1]{Antti Hannukainen}
\author[2]{Nuutti Hyv\"onen}
\author[3]{Vigdis Toresen}
\affil[1]{Department of Mathematics and Systems Analysis, Aalto University, P.O. Box 11100, 00076 Helsinki, Finland.}
\affil[2]{Department of Mathematics and Systems Analysis, Aalto University, P.O. Box 11100, 00076 Helsinki, Finland.}
\affil[3]{Department of Mathematics and Systems Analysis, Aalto University, P.O. Box 11100, 00076 Helsinki, Finland.}


\begin{abstract}
This work considers using reduced basis techniques in connection to (smoothened) total variation regularization in electrical impedance tomography, but analogous ideas can also be used for other inverse elliptic boundary value problems. It is demonstrated that resorting to reduced bases can speed up a reconstruction algorithm based on combining the lagged diffusivity algorithm with sequential linearizations and preconditioned LSQR iteration without any significant loss of reconstruction quality or of the edge-enhancing nature of total variation regularization. The ideas are numerically tested in three dimensions on unstructured finite element meshes with both simulated and experimental data, resulting in online reconstruction times of only a few seconds on a standard laptop computer. 
\end{abstract}

\begin{document}

\flushbottom
\maketitle
\thispagestyle{empty}

\section{Introduction}
\textit{Electrical impedance tomography} (EIT) is a non-invasive imaging technique in which information about the internal conductivity of an object is reconstructed based on injecting currents through electrodes placed on the surface and measuring the resulting electrode voltages. We model EIT with the \textit{complete electrode model} (CEM) that accounts for the shapes of the electrodes and the resistive layers at the electrode-object interfaces. The CEM is capable of predicting the operation of real-world EIT devices up to measurement precision \cite{CEM,ExistenceUniqueness}. For general information on EIT and its mathematical idealization, i.e., the Calder\'on problem, we refer to the review articles \cite{Borcea02,Uhlmann09} and the references therein. 

This work tackles EIT in the framework of Bayesian inversion under the prior information that the imaged object is composed of embedded inclusions in a constant background. The prior knowledge is encoded in a smoothened {\em total variation} (TV) type prior. The main goals are (i) to present an accelerated MATLAB-based variant of the edge-enhancing reconstruction algorithm introduced in \cite{Edge-Enhancing}, allowing reconstructions on three-dimensional unstructured {\em finite element} (FE) meshes in a couple of tens of seconds on a standard laptop computer, and (ii) to numerically demonstrate that reduced basis techniques can be combined with a TV prior or regularization to (further) speed up the online phase of the algorithm without compromising reconstruction quality or the capability to reconstruct jumps in the conductivity. Although we only consider EIT, the presented ideas could also be applied to other inverse elliptic (boundary value) problems, cf.~\cite{Hannukainen16a,Hannukainen16b}. It should also be noted that there exist a plethora of other works applying TV regularization to EIT; see, e.g., \cite{Borsic10,Chung05,Kaipio00,Kolehmainen98}.

The studied reconstruction algorithm aims to find a \textit{maximum a posteriori} (MAP) estimate for the discretized conductivity, which is a task that can be formulated as the minimization of a Tikhonov functional that is non-quadratic in both of its terms. The nonquadraticity is due to the nonlinearity of the studied inverse problem and the use of a TV prior, respectively. The algorithm is composed of two nested iterations: the outer loop corresponds to sequential linearizations of the measurement model, whereas the inner loop applies LSQR \cite{LSQR,Paige82a,Paige82b} iteration and prior-conditioning \cite{Calvetti07,Calvetti12,Calvetti05} to taking one lagged diffusivity step \cite{lagged} for the linearized TV regularized problem. As a deviation from the material in \cite{Edge-Enhancing}, we resort to a recent idea introduced in~\cite{Jaaskelainen25} and do not reconstruct the contact resistances but project the to-be-inverted equations onto the orthogonal complement of the range of the Jacobian matrix of the electrode measurements with respect to the contacts; see also~\cite{Calvetti25}. Such a simplified approach has been demonstrated to produce good-quality reconstructions of the internal conductivity (almost) independently of the accuracy of the initial guess for the conductivity and contact resistances at which the employed projection is computed~\cite{Jaaskelainen25}. Compared to \cite{Edge-Enhancing}, the main speed-up in the performance of the basic version of the algorithm is due to (i) precomputing quantities that do not depend on the conductivity in the offline phase of the algorithm, i.e., when the imaging geometry is already known but the measurements are not yet in hand, and (ii) paying careful attention to optimal employment of data structures in MATLAB.

The reduced basis version of the algorithm is based on building a {\em proper orthogonal decomposition} (POD) \cite{POD} for the solutions of the CEM forward problem, or more precisely, only for the electric potentials induced {\em inside} the body by the applied electrode currents. The POD requires acquiring a library of snapshots of the potentials; these are computed by solving the forward problem for a number of conductivities and contact resistance drawn from log-normal distributions during the offline phase of the algorithm. In particular, forming the POD does not utilize the prior information on the piecewise constant structure of the imaged object. The POD allows replacing the large sparse FE system matrices in the reconstruction algorithm by much smaller but dense ones obtained via projecting onto the orthonormal basis produced by the POD. This considerably accelerates computing both the forward solutions and the utilized Jacobians with respect to the discretized conductivity. However, it is essential to notice that no reduced basis is used for evaluating the TV prior, that is, in terms of regularizaton, the POD is used for the data fit term but not for the penalty term in the Tikhonov functional. POD has been used in connection to EIT,~e.g.,~in~\cite{reduced_eit,EIT_POD_1, EIT_POD_2}, but we are not aware of any previous works that have successfully combined any reduced basis technique with TV regularization for EIT or other inverse elliptic boundary value problems.

This text is organized as follows. Section~\ref{sec:CEM} recalls the CEM, and Section~\ref{sec:RB} explains how a reduced basis can be built for the associated forward solutions. The adopted Bayesian inversion framework with a smoothened TV prior is introduced in Section~\ref{sec:Bayes}, and the projected lagged diffusivity algorithm for seeking a MAP-like solution for the considered reconstruction problem is described in Section~\ref{sec:ProjectedLaggedDiffAlg}. The implementation of the reconstruction algorithm, both with and without reduced basis, is considered in Section~\ref{sec:implementation}. Sections~\ref{sec:numerics} and~\ref{sec:conclusion} present, respectively, the numerical experiments and the concluding remarks.

\section{Complete electrode model}
\label{sec:CEM}

CEM is a commonly used mathematical model for EIT that has been shown to give results in agreement with experiments \cite{CEM}. Although EIT measurements are usually performed using alternating current, which leads to a model with complex-valued current patterns and electric potentials~\cite{Borcea02}, we ignore the phase information and use the amplitudes of electrode currents and voltages as inputs and measurements, respectively, for a forward model with real-valued conductivity and contact resistances. Such an approximation/simplification has been successfully used in numerous works on experimental EIT. To allow straightforward comparison with the material in \cite{Edge-Enhancing}, we adopt the basic version of CEM that assumes the contact resistivity is distributed uniformly over the surface of an individual electrode, even though the employed idea to project away the contribution of contacts was introduced in \cite{Jaaskelainen25} for the smoothened version of the CEM that allows a varying surface resistivity over the object boundary~\cite{Hyvonen17b}.

Let $\Omega \subset \R^3$ be a bounded domain with a sufficiently smooth boundary $\partial \Omega$, which is partially covered with $M\in \mathbb{N}\setminus \{1\}$ electrodes represented as non-overlapping open connected subsets $E_m \subset \partial \Omega$, $m=1,\hdots, M$. Their union is denoted by $E$. The conductivity distribution $\sigma \in L^\infty (\Omega)$ is assumed to satisfy $\sigma \geq c > 0$ almost everywhere in $\Omega$. As mentioned above, the contact resistances $z\in \R^M$ are modeled as constant on each electrode and assumed to be bounded away from zero $z_m \geq c > 0$, $m=1,\hdots ,M$. 

In an EIT measurement, currents $I_m$, $m=1, \dots, M$, are injected through the electrodes and the resulting voltages $U_m$, $m=1, \dots, M$, on the electrodes are measured. Under the assumption that the domain $\Omega$ contains no sources or sinks, conservation of charge dictates that the current pattern belongs to the space
\begin{equation}\label{Rdiamond}
    \R_{\diamond}^M= \bigg\{ y \in \R^M \ \Big| \ \sum_{i=1}^M y_i = 0 \bigg\}.
\end{equation}
The ground potential is chosen so that the vector of voltages measured on the electrodes also satisfies $U\in \R_{\diamond}^M$. Together with the electric potential inside the body $u$, these can be determined as the
weak solution $(u,U)\in H^1(\Omega)\oplus \R_{\diamond}^M$ to the elliptic boundary value problem
\begin{equation}
\label{eq:cemeqs}
\begin{array}{rll}
    \nabla \cdot (\sigma \nabla u) \!\! \!\! &= 0 \qquad &\text{in } \Omega, \\[2mm]
     {\displaystyle\frac{\partial u}{\partial \nu}} \!\! \!\! &= 0 \qquad &\text{on } \partial\Omega \setminus E,\\[3mm]
    {\displaystyle u + z_m\sigma \frac{\partial u}{\partial \nu}} \!\! \!\! &= U_m \quad &\text{on } E_m,\quad m=1,\hdots ,M,\\[3mm]
    {\displaystyle \int_{E_m} \sigma \frac{\partial u}{\partial \nu}} \, {\rm d} S \!\! \!\! &= I_m, \quad & m=1, \hdots , M.
    \end{array}
\end{equation}
Here, $\nu$ is the unit exterior normal vector to $\partial \Omega$. A proof of the existence and uniqueness of the solution to \eqref{eq:cemeqs} was given in \cite{ExistenceUniqueness}.

We solve \eqref{eq:cemeqs} numerically using the {\em finite element method} (FEM) with piecewise linear (first order) basis functions $\{\psi_j\}_{j=1}^N$ on a 3D tetrahedral mesh with appropriate refinements close to the electrode boundaries. By abuse of notation, we denote by $(u,U) \in {\rm span} \{\psi_j\}_{j=1}^N \oplus {\rm span} \{\mathbf{c}^{(m)}\}_{m=1}^{M-1}$ this approximate FE solution, with $\{\mathbf{c}^{(m)}\}_{m=1}^{M-1}$ being the employed orthonormal basis for $\mathbb{R}_\diamond^M$. We also use a piecewise linear FE basis $\{\varphi_j\}_{j=1}^n$ to represent the conductivity~$\sigma$. Although in principle the two bases, $\{\psi_j\}_{j=1}^N$ and $\{\varphi_j\}_{j=1}^n$, can be different, we use the same one for $\sigma$ and $u$ in the numerical experiments. However, since we only use the reduced basis on the representation of the internal potential, it is useful to have different symbols for the bases. In what follows, any function presented in the basis $\{\psi_j\}_{j=1}^N$ or $\{\varphi_j\}_{j=1}^n$ is identified with the corresponding vector of nodal values, or degrees of freedom, on the associated FE mesh. As an example, $u$ can be interpreted either as a element in $H^1(\Omega)$ or as a vector in $\R^N$. It should be clear from the context which is meant. Details on how the FE system matrix $A\in \R^{(N+M-1)\times (N+M-1)}$ for \eqref{eq:cemeqs} is constructed can be found in \cite{Vauhkonen_thesis} or \cite{reduced_eit}. Note that the degrees of freedom in the FE solution are ordered so that the first $N$ correspond to $u$ and the last $M-1$ to $U$.

The inverse problem of EIT is to reconstruct $\sigma$ from noisy voltage measurements on the electrodes. Measurements are typically performed for $M-1$ linearly independent current patterns. Not accounting for the measurement noise, such a number is necessary and sufficient for collecting all information available via EIT measurements with $M$ electrodes due to the linear dependence of the solution pair $(u,U)$ on $I$. However, as some of our numerical tests do not use all electrodes for current injection and, on the other hand, utilizing more than $M-1$ current patterns can be viewed as a technique for reducing measurement noise, we let $L \in \mathbb{N}$ denote the number of current patterns. The applied current patterns and electrode voltages predicted by the FE discretized CEM are collected into the matrices
\begin{equation}
    \mathcal{I}= 
    \begin{bmatrix}I^{(1)}, \dots ,   I^{(L)}\end{bmatrix} \in \R^{M\times L}.
\end{equation}
and 
\begin{equation}\label{Ucal}
    \mathcal{U}(\sigma, z; I^{(1)},\hdots ,I^{(L)}) = \begin{bmatrix} U(\sigma, z; I^{(1)}), \hdots , U(\sigma, z; I^{(L)}) \end{bmatrix} \in \R^{M \times L},
\end{equation}
respectively. The linear operation of column-wise vectorization of a matrix is denoted by ${\it vec}$, and we define
\begin{equation}
\label{eq:U_to_U}
\mathbf{U}(\sigma, z) = {\it vec} ( \mathcal{U}(\sigma, z)) \in \R^{LM} 
\end{equation}
to be the vectorized version of the predicted measurements. In \eqref{eq:U_to_U}, the dependence on the current patterns is suppressed, which is the convention used in the rest of this paper.

Note that the electrode potential measurements also depend on the contact resistances~$z$ that are nuisance parameters: their values are not of primary interest, but ignoring them in the inversion process may lead to severe artifacts especially if absolute EIT measurements are considered, as is the case in this work; see,~e.g.,~\cite{Heikkinen02,Vilhunen02}. The standard approach would be to include $z$ as a second unknown for the inverse problem, but here we instead resort to a technique for approximately projecting away the contribution of the contact resistances~\cite{Jaaskelainen25}.


\section{Reduced basis}\label{sec:RB}

This section briefly reviews how reduced basis can be utilized when numerically solving~\eqref{eq:cemeqs}. For a more comprehensive exposition, see \cite{reduced_eit}.

\subsection{Review of subspace surrogate methods}

Let $A\in \R^{m\times m}$ be a symmetric positive definite matrix, $\mathbf{b}\in \R^{m}$ and $\mathbf{x}\in \R^{m}$ the solution to
\begin{equation}\label{Axb}
    A \mathbf{x} = \mathbf{b}.
\end{equation}
Let $\mathcal{W}\subset \R^m$ be a subspace with a basis $\{\mathbf{q}_i\}_{i=1}^k$ with $k<m$, and introduce the associated basis matrix $Q=\begin{bmatrix}
    \mathbf{q}_1,\hdots , \mathbf{q}_k\end{bmatrix}$.
Instead of (\ref{Axb}), we may consider solving the smaller surrogate system
\begin{equation}
    Q^\top \! A Q \hat{\mathbf{x}} = Q^\top\mathbf{b},
\end{equation}
and approximate $\mathbf{x}$ with $\tilde{\mathbf{x}} = Q \hat{\mathbf{x}} $. If $A$ originates from a FE discretization, as in our numerical studies, it is sparse but typically of a very large dimension. On the other hand, the surrogate matrix $Q^\top \! A Q $ is usually dense, but much smaller if $k \ll m$.

Let us then consider the special case of interest that $A \in \R^{(N + M-1) \times (N + M-1)}$ is the FE stiffness matrix for the CEM introduced in Section~\ref{sec:CEM}. Typically, the dimension $N$ of the FE subspace for $H^1(\Omega)$ is large, while the number of electrodes $M$ is small, e.g., $M=16$ or $M=32$. It therefore makes sense to only reduce the basis for the inner potentials and not the electrode potentials. Thus, we deviate from the above introduced structure for the matrix $Q$ and write it instead in a block form as
\begin{equation}
\label{eq:tildeQ}
    Q=\begin{bmatrix}
        \widehat{Q} & 0 \\[0.5mm] 0 & \mathrm{I}
    \end{bmatrix},
\end{equation}
where $\widehat{Q}\in \R^{N\times k}$ is the basis matrix used to reduce the inner potentials and $\mathrm{I}\in \R^{(M-1)\times (M-1)}$ is an identity matrix applied to the electrode potentials.

\subsection{Basis generation}
Building the reduced basis matrix $Q$ is part of the offline stage and not considered time-critical. In other words, it is assumed that there is sufficient information on the imaging setup well in advance to allow building the reduced basis prior to performing the measurements. It was numerically shown in \cite{reduced_eit} that different methods for creating the reduced subspace give similar reconstructions for the EIT inverse problem.

\subsubsection{Proper orthogonal decomposition}

We build the reduced basis matrix $\widehat{Q}$ in \eqref{eq:tildeQ} using POD; see \cite{POD} for an introduction. POD has previously been applied to the EIT inverse problem in \cite{reduced_eit,EIT_POD_1, EIT_POD_2}.

The forward model is solved $K > k$ times for randomly drawn $\sigma$ and $z$ --- in our numerical experiments, from log-normal distributions as explained in the next subsection. The resulting snapshots are collected in the matrix $Y_K = [u_1(\sigma_1, z_1), \hdots, u_K(\sigma_k, z_K)] \in \R^{N \times K}$. An orthonormal basis for the $k$-dimensional POD subspace, i.e.~the columns of $\widehat{Q}$, are constructed as eigenvectors corresponding to the $k$ largest eigenvalues of the \textit{correlation matrix} $Y_K Y_K^\top \in \R^{N \times N}$. Note that the correlation matrix does not need to be formed explicitly, but one can exploit a singular value decomposition for $Y_K$. The singular values have been observed \cite{reduced_eit} to decay exponentially, enabling the use of randomized linear algebra to efficiently compute the range of the needed truncated singular value decomposition. More specifically, we use the randomized range finder~\cite[Algorithms~4.1 \& 5.1]{Randomized}.

\subsubsection{Sampling}
During training, the samples for $\sigma$ and $z$ are drawn from log-normal distributions, the specifications of which are given in Section~\ref{sec:numerics}. The reason for resorting to log-normal rather than normal distributions is the required non-negativity of $\sigma$ and $z$. In particular, the teaching sample for $\sigma$ has very different characteristic compared to the sharply defined inclusions we aim to reconstruct with the TV prior. That is, we do not include the prior information on the piecewise constant structure of $\sigma$ in the reduced basis matrix~$Q$ for the interior potential. However, we assume that we have a reasonable estimate for the background conductivity $\sigma_0$, the logarithm of which is used as the mean for the underlying normal distribution. In addition, we assume to have some kind of an understanding on the characteristic distances of structures inside $\Omega$ in order to choose a correlation length for the employed log-normal density for $\sigma$. However, our numerical examples demonstrate that this choice is not crucial for the functionality of the reduced basis approach.


\section{Bayesian framework and prior}
\label{sec:Bayes}

The voltage measurements on the electrodes are modeled as
\begin{equation}
\label{eq:meas_model}
    \mathbf{V} = \mathbf{U}(\sigma, z) + \eta,
\end{equation}
where the noise vector $\eta \in \R ^{LM}$ is a realization of a Gaussian random variable with zero mean and a known positive definite covariance matrix $\Gamma \in \R ^{LM\times LM}$. The likelihood of the measurement $\mathbf{V}$ given the parameters $\sigma$ and $z$ is thus
\begin{equation}\label{likelihood}
    p(\mathbf{V} \mid \sigma, z) \propto \exp \! \Big(-\frac{1}{2}(\mathbf{V}-\mathbf{U}(\sigma, z))^\top \Gamma^{-1}(\mathbf{V}-\mathbf{U}(\sigma, z) ) \Big).
\end{equation}
We use the symbol $p$ to refer to a general (unnormalized) probability density.

The prior information on the conductivity is that it consists of an approximately homogeneous background  with sharply defined inclusions. This is encoded into the prior probability density assigned to the conductivity,
\begin{equation}
    p(\sigma) \propto \exp (-\alpha R(\sigma)),
\end{equation}
where $\alpha > 0$,
\begin{equation}
    R(\sigma) = \int_\Omega r (\lVert \nabla \sigma (x)\rVert) \, {\rm d} x,
\end{equation}
and $\lVert \, \cdot \, \rVert$ denotes the Euclidean norm. We use a smoothened TV prior, meaning that
\begin{equation}
    r(t) = \sqrt{T^2 + t^2} \approx |t|.
\end{equation}
The small parameter $T>0$ ensures the differentiability of $r$. This type of prior promotes sharp changes or edges over slow changes or oscillations of the conductivity $\sigma$. For the contact resistances $z$, we use an uninformative prior that gives equal prior probability to all realizations of $z$.

Bayes' rule gives the posterior probability density
\begin{equation}
\label{eq:posterior}
     p(\sigma, z \mid \mathbf{V}) \propto  p(\mathbf{V} \mid \sigma, z) \, p(\sigma),
\end{equation}
where the hidden multiplicative constant does not depend on $\sigma$ or $z$. Determining the (possibly non-unique) MAP estimate for this posterior is equivalent to minimizing the Tikhonov functional
\begin{equation}\label{Phi}
    \Phi (\sigma, z) := \frac{1}{2}(\mathbf{V}-\mathbf{U}(\sigma, z))^\top \Gamma^{-1}(\mathbf{V}-\mathbf{U}(\sigma, z) ) + \alpha R(\sigma ) 
\end{equation}
that is non-quadratic in both the term coming from the likelihood \eqref{likelihood} and the penalty term originating from the prior.

\begin{remark}
As the nodal values of the conductivity $\sigma \in \R_+^n$ as well as the contact resistances $z \in \R_+^M$ are positive, the respective (improper) prior densities should include positivity priors,~i.e.,~multiplicative terms that take the value $1$ only if all components of the considered vector are positive and the value $0$ otherwise. For simplicity, these terms have been omitted from the above analysis. The omission of the positivity prior for the contact resistances does not affect the reconstruction algorithm introduced in the following section since their contribution to the forward model is approximately suppressed via a projection \cite{Jaaskelainen25}. On the other hand, the positivity prior for the nodal conductivities is implicitly accounted for in the reconstruction algorithm by replacing possible negative components by a small positive number after each iteration round. 
\end{remark}

\section{Projected lagged diffusivity algorithm}\label{sec:ProjectedLaggedDiffAlg}

Lagged diffusivity iteration is a ﬁxed-point algorithm introduced in \cite{lagged} for TV denoising in image processing. However, it can also be generalized for other finite-dimensional inverse problems \cite{LSQR,Ascher06}, including nonlinear problems via sequential linearization \cite{Hannukainen16a, Hannukainen16b,Edge-Enhancing}. The algorithm version presented here is essentially the one from \cite{Edge-Enhancing} with the difference that we project out the dependence on $z$ at the start and make no further attempt at reconstructing the contact resistances. When comparing our results to those of \cite{Edge-Enhancing}, we find that the projection does not much affect the quality of the internal conductivity reconstruction, but it somewhat simplifies the algorithm presented here.

\subsection{Sequential Linearization}

Without employing the projection with respect to the contact resistances, the iterative algorithm is based on applying (sequential) linearization of the forward model to \eqref{Phi}, which results in the task of minimizing
\begin{equation}\label{Phi_i}
    \Phi^{(i)} (\sigma, z) := \frac{1}{2}\big(y^{(i)}-J_\sigma^{(i)} {\sigma} -J_z^{(i)}z \big)^\top \Gamma^{-1} \big(y^{(i)}-J_\sigma^{(i)} {\sigma} -J_z^{(i)}z \big) + \alpha R(\sigma ),
\end{equation}
where $(\sigma^{(i)}, z^{(i)})$ is the current estimate for the conductivity and contact resistances, 
\begin{equation}\label{yi}
y^{(i)}=\mathbf{V}-\mathbf{U}({\sigma}^{(i)}, z^{(i)})+J_\sigma^{(i)} {\sigma}^{(i)} + J_z^{(i)}z^{(i)},
\end{equation} 
and $J_\sigma^{(i)}$ and $J_z^{(i)}$ are the Jacobian matrices of the map $({\sigma}, z) \mapsto \mathbf{U}({\sigma}, z)$ with respect to ${\sigma}$ and $z$, respectively, evaluated at $({\sigma}^{(i)}, z^{(i)})$. The linearization is called sequential because the Jacobians and the basis point for the linearization are updated at each step. In the next section, we will explain why it is only necessary to evaluate $J_z$ once and use it to project out the dependence on $z$.

We assume that we have a good initial guess for the background conductivity $\sigma_0 \in \R^n$ and a homogeneous estimate $z_0 = \zeta_0 \begin{bmatrix} 1, \dots, 1 \end{bmatrix}^\top \in \R_+^M$ for the contact resistances. If such are not available, they can be estimated at the start of the online stage \cite{Edge-Enhancing}, which would somewhat add to the computation time. Since this is not the focus of the present work, we assume that $\sigma_0$ and $z_0$ are available at the start. Note that the homogeneous estimate $z_0$ for the contact resistances need not be accurate for the algorithm to function;~cf.~\cite{Jaaskelainen25}.


\subsection{Projected equations}
\label{sec:projection}

The orthogonal projection onto the orthogonal complement of the range $\mathcal{R}(J_z)$ is
\begin{equation}\label{Pz}
    P = \mathrm{I} - J_z (J_z^\top \! J_z)^{-1} J_z^\top,
\end{equation}
where $ \mathrm{I}\in \R^{LM \times LM}$ is the identity matrix of the appropriate size. The numerical considerations in \cite{Jaaskelainen25} for similar projection operators demonstrate that the range of $P$ is almost independent of $\sigma$ and $z$. This suggests an approach to EIT reconstruction where $P$ is formed only once using the initial guess $(\sigma_0,z_0)$, and the corresponding $P=P(\sigma_0, z_0)$ is then used to approximately project out the dependence on $z$ to allow only reconstructing the internal conductivity. The formation of $P$ can be part of the offline stage, assuming the initial guesses are available; note also that constructing $P$ is computationally cheap, if $J_z$ is available, as the number of electrodes $M$ is small.


The underlying idea in using $P$ is to replace the measurement model \eqref{eq:meas_model} by a projected version
\begin{equation}
\label{eq:proj_meas_model}
    P  \mathbf{V} = P \mathbf{U}(\sigma, z) + P \eta,
\end{equation}
which corresponds to only considering the available data on the orthogonal complement of $\mathcal{R}(J_z)$. The noise $P\eta$ is still zero-mean Gaussian but supported on $\mathcal{R}(J_z)^\perp$ with the (only positive semidefinite) covariance matrix $P \Gamma P^\top = P \Gamma P$, as $P^\top = P$ for an orthogonal projection. Moreover, restricted on $\mathcal{R}(J_z)^\perp$, the inverse of the covariance matrix of $P\eta$ is $(P \Gamma P)^\dagger = P \Gamma^{-1} P$ because the pseudoinverse $P^\dagger$ is $P$ itself.

Exactly as in the non-projected case (see \eqref{eq:posterior}, \eqref{Phi} and \eqref{Phi_i}), it straightforwardly follows that finding a MAP estimate for the posterior $p(\sigma, z \mid P \mathbf{V})$ corresponds to minimizing 
\begin{equation}\label{eq:PPhi}
    \Phi_P(\sigma, z) := \frac{1}{2}\big(P\mathbf{V}-P\mathbf{U}(\sigma, z)\big)^\top P\Gamma^{-1}P\big(P\mathbf{V}-P\mathbf{U}(\sigma, z) \big) + \alpha R(\sigma ).
\end{equation}
According to numerical tests in \cite{Jaaskelainen25},  $\mathcal{R} (J_z (\sigma_0, z_0)) \approx \mathcal{R} (J_z (\sigma, z))$ for a wide variety of $(\sigma, z)$, which indicates $P J_z (\sigma, z) \approx 0$. Trusting this approximation, the linearization of the forward model around $(\sigma^{(i)}, z^{(i)})$ in \eqref{eq:PPhi} yields
\begin{align} 
\label{eq:PPhi_i:init}
\Phi_P(\sigma, z) &\approx \frac{1}{2}\big(P y^{(i)}-P J_\sigma^{(i)} {\sigma} - P J_z^{(i)}z \big)^\top P \Gamma^{-1} P \big(P y^{(i)}- P J_\sigma^{(i)} {\sigma} - P J_z^{(i)}z \big) + \alpha R(\sigma ) \nonumber \\[1mm]
&\approx \frac{1}{2}\big(P y^{(i)}- P J_\sigma^{(i)} {\sigma} \big)^\top P \Gamma^{-1} P \big(P y^{(i)}- PJ_\sigma^{(i)} {\sigma} \big) + \alpha R(\sigma ).
\end{align}
Using the approximation 
\begin{align}
\label{eq:b_vector}
P y^{(i)} &= P \mathbf{V} - P \mathbf{U}({\sigma}^{(i)}, z^{(i)})+ P J_\sigma^{(i)} {\sigma}^{(i)} + P J_z^{(i)}z^{(i)} \nonumber \\[1mm]
& \approx
 P \big (\mathbf{V} -  \mathbf{U}({\sigma}^{(i)}, z^{(i)})+  J_\sigma^{(i)} {\sigma}^{(i)} \big) =: P y^{(i)}_P
\end{align}
and the identities $P = P^\top$ and $PP = P$
on the right-hand side of \eqref{eq:PPhi_i:init}, we finally arrive at the approximate Tikhonov functional
\begin{equation}\label{PPhi_i}
   \Phi_P^{(i)} (\sigma) := \frac{1}{2} \big (y^{(i)}_P-J_\sigma^{(i)} {\sigma} \big )^\top P\Gamma^{-1}P \big(y^{(i)}_P-J_\sigma^{(i)} {\sigma} \big) + \alpha R(\sigma ).
\end{equation}
Since our projected iterative reconstruction algorithm finds $(\sigma^{(i+1)}, z^{(i+1)})$ by minimizing $\Phi_P^{(i)}(\sigma)$ with respect to $(\sigma, z)$, but $\Phi_P^{(i)}(\sigma)$ no longer depends on $z$, we deduce that $z^{(i)} = z^{(0)} = z_0$ remains fixed during the iteration. Hence, $y^{(i)}_P$ and $J_\sigma^{(i)}$ can be redefined in \eqref{PPhi_i} as
\begin{equation}
\label{eq:J0_and_yP0}
J_\sigma^{(i)} = J_\sigma(\sigma^{(i)}, z_0) \qquad \text{and} \qquad y^{(i)}_P = \mathbf{V} -  \mathbf{U}({\sigma}^{(i)}, z_0) +  J_\sigma^{(i)} {\sigma}^{(i)},
\end{equation}
respectively.

The changes in $\Phi_P^{(i)}$ compared to the non-projected version $\Phi^{(i)}$ defined in \eqref{PPhi_i} are that the term $J_z^{(i)}$ is dropped from the quadratic form and from the definition of the target vector (cf.~\eqref{yi} and \eqref{eq:J0_and_yP0}), and the inverse covariance matrix of the noise $\Gamma^{-1}$ is replaced by its projected version $P\Gamma^{-1}P$. For more detailed Bayesian analysis of the adopted projection idea in the framework of finite-dimensional linear inverse problems, consult~\cite{Calvetti25}.

\subsection{MAP estimate}

The necessary condition for a local minimum of (\ref{PPhi_i}) is $\nabla \Phi_P^{(i)} =0$, which gives
\begin{equation}\label{nablaPhi_i}
    (J_\sigma^{(i)})^\top P\Gamma^{-1}P J_\sigma^{(i)}\sigma + \alpha (\nabla R)(\sigma ) = (J_\sigma^{(i)})^\top P\Gamma^{-1}P y^{(i)}_P.
\end{equation}
This is, in fact, also a sufficient condition due to the convexity of $\Phi_P^{(i)}$ if the left-hand side of  \eqref{nablaPhi_i} does not vanish for a constant $\sigma$, which can be guaranteed, e.g., by imposing a suitable Dirichlet boundary condition for $\sigma$; see,~e.g.,~\cite{Pohjavirta21}.  The gradient term in \eqref{nablaPhi_i} can be represented in a matrix form as \cite{LSQR}
\begin{equation}
    (\nabla R)(\sigma ) = H(\sigma ) \sigma.
\end{equation}
The entries of the matrix $H\in \R^{n\times n}$ are given by
\begin{equation}\label{H}
    H_{ij}(\sigma) = \int_{\Omega}  \frac{\nabla\varphi_i(x) \cdot \nabla \varphi_j(x)}{\sqrt{T^2 + \| \nabla \sigma (x)\|^2}} \, {\rm d} x, \quad i,j=1,\hdots, n,
\end{equation}
where, as introduced in Section~\ref{sec:CEM}, $\{\varphi_i\}_{i=1}^n$ is the FE basis used to discretize the conductivity such that
\begin{equation}
    \sigma = \sum_{j=1}^n \sigma_j \varphi_j.
\end{equation}
$H$ can be interpreted as the FE stiffness matrix of an elliptic partial differential operator
\begin{equation}
\label{eq:c}
    - \nabla \cdot (c_\sigma (x) \nabla ( \, \cdot \,)),
\end{equation}
where
\begin{equation}
    c_\sigma : x \mapsto  \frac{1}{\sqrt{T^2 + \| \nabla \sigma (x)\|^2}}, \quad \Omega \rightarrow \R_+,
\end{equation}
accompanied by a homogeneous Neumann boundary condition on $\partial \Omega$. 

Since adding a constant term in the target function would not affect the result after applying the operator \eqref{eq:c} and the natural boundary condition, $H$ cannot be invertible as such --- in fact, it follows from standard FE analysis that $H$ is positive semidefinite with its nullspace spanned by the constant vector. To remedy the problems this would cause in subsequent analysis, we impose artificial Dirichlet boundary condition on a subset of $\partial \Omega$ to ensure the positive definiteness of $H$ and thus the existence of $H^{-1}$. We will return to how this is exactly done below; at this stage, we simply note that a Dirichlet boundary condition fixes a certain number of degrees of freedom in FE representations, and thus in the following, we abuse the notation by interpreting $H \in \R^{n' \times n'}$, $J_\sigma \in \R^{LM \times n'}$, and $\sigma \in \R^{n'}$ for a certain natural number $n' < n$. An alternative way to ensure invertibility of $H$ would be to add a small positive constant to the diagonal elements of $H$. This is used,~e.g.,~in \cite{Jaaskelainen25}.


\subsection{Lagged diffusivity}
\label{sec:LD}

The equation to be solved is now in a matrix form,
\begin{equation}\label{nonlineq}
\big((J_\sigma^{(i)})^\top P\Gamma^{-1}PJ_\sigma^{(i)} + \alpha H(\sigma^{(i+1)} ) \big) \sigma^{(i+1)}= (J_\sigma^{(i)})^\top P\Gamma^{-1}Py^{(i)}_P,
\end{equation}
which is still a nonlinear equation since $H(\sigma^{(i+1)})$ depends on $\sigma^{(i+1)}$. To solve this, we use the lagged diffusivity iteration. The idea is to iteratively solve \eqref{nonlineq} by ``lagging'' the nonlinear term,~i.e.,~using its evaluation at the previous step of the iteration. This leads, in principle, to a nested iteration within the sequential linearization:

\begin{enumerate}
    \item Set $\sigma^{(i+1)}_{(0)} = \sigma^{(i)}$ and $j=0$.
    \item Iteratively solve $\sigma^{(i+1)}_{(j+1)}$ from the linear equation
\begin{equation}
\big ((J_\sigma^{(i)})^\top P\Gamma^{-1}PJ_\sigma^{(i)} + \alpha H(\sigma^{(i+1)}_{(j)} ) \big) \sigma^{(i+1)}_{(j+1)}= (J_\sigma^{(i)})^\top P\Gamma^{-1}P y^{(i)}_P,
\end{equation}
which is uniquely solvable due to the positive definiteness of $H(\sigma^{(i+1)}_{(j)})$.
\item If the stopping criterion is satisfied, assign $\sigma^{(i+1)}=\sigma^{(i+1)}_{(j+1)}$. Otherwise, increment $j\rightarrow j+1$ and go back to step 2.
\end{enumerate}
For analysis on the convergence of the lagged diffusivity iteration toward the minimizer of \eqref{PPhi_i}, we refer to \cite{Chan99,Dobson97,Pohjavirta21}. 

\subsection{Priorconditioning}

We now look at a single step of the lagged diffusivity iteration nested inside the sequential linearizations. We drop the indices for readability, that is, we denote
$\sigma = \sigma^{(i+1)}_{(j+1)}$ and $H=H(\sigma^{(i+1)}_{(j)} )$. To further simplify the notation, we introduce $B\in \R^{LM\times n'}$ and $b\in \R^{LM}$ defined by
\begin{equation}\label{Bb}
    B = \Gamma^{-1/2}PJ_\sigma^{(i)} \qquad \text{and} \qquad b = \Gamma^{-1/2}Py_P^{(i)},
\end{equation}
where $\Gamma^{-1/2}$ is a Cholesky factor of $\Gamma^{-1}$. In the numerical experiments, $\Gamma$ is diagonal, which justifies the inverse square root notation. In this notation, we have recognizably a regularized normal equation,
\begin{equation}\label{normal_equation}
    (B^\top \! B + \alpha H) \sigma = B^\top b,
\end{equation}
which is equivalent to the minimization of the least squares functional
\begin{equation}
\label{eq:LQ}
    \lVert B\sigma - b\rVert^2 + \alpha \sigma^\top \! H \sigma,
\end{equation}
where $H$ defines the form of penalization. We solve (a version of) the normal equation~\eqref{normal_equation} by combining prior conditioning, LSQR and the Morozov discrepancy principle. An alternative would be to exploit the thin structure of $B$ via the Woodbury matrix identity as in \cite{Jaaskelainen25}.

There are some problems we need to address before attempting to solve \eqref{normal_equation}. First, there is the free parameter $\alpha$ that would need to be determined based on prior knowledge or eliminated in some other way. Second, the matrix $B^\top \! B$ is more ill-conditioned than $B$. To tackle these, we precondition by multiplying \eqref{normal_equation} from the left by $H^{-1} = (H^{-1})^\top$, which results in
\begin{equation}
    \big((H^{-1})^\top \! B^\top \! B + \alpha \mathrm{I} \big)  {\sigma} = (H^{-1})^\top \! B^\top {b}.
\end{equation}
 Preconditioning with the inverse prior covariance matrix, which is a possible Bayesian interpretation for $H$ in light of \eqref{eq:LQ}, was dubbed priorconditioning in \cite{Calvetti05}. We assume that this priorconditioning encodes the essential available prior information on $\sigma$ (cf.~\eqref{eq:LQ}) and set $\alpha = 0$ in the following. We also make the change of variables $\tilde{\sigma} = \sigma - \sigma_0$ and note that the Jacobian with respect to $\tilde{\sigma}$ is still $J_\sigma$, so that $B$ defined in \eqref{Bb} remains unchanged. However, the target vector $b$ needs to be shifted, that is, 
\begin{equation}\label{btilde}
  \tilde{b} = b - B\sigma_0 = \Gamma^{-1/2}P\tilde{y}_P^{(i)}, \qquad \text{with} \qquad  \tilde{y}_P^{(i)}=\mathbf{V}-\mathbf{U}(\sigma^{(i)}, z_0) +J_\sigma^{(i)} \tilde{\sigma}^{(i)}.
\end{equation}

These changes leave us with the equation
\begin{equation}\label{preconditioned}
    (H^{-1})^\top \! B^\top \! B \tilde{\sigma} = (H^{-1})^\top \! B^\top \tilde{b},
\end{equation}
which can be efficiently solved by applying a preconditioned version \cite{LSQR} of LSQR \cite{Paige82a,Paige82b} in the same way as it was used in \cite{Edge-Enhancing}. LSQR uses $B$ instead of the more ill-conditioned $B^\top \! B$ but is equivalent to the conjugate gradient method for normal equations. The LSQR algorithm only requires the ability to operate with the inverse of $H$ on a given vector; that is, there is no need to introduce some specific factorization for $H$ (or to explicitly form $H^{-1}$ for that matter). Be that as it may, the LSQR algorithm from \cite{LSQR} is equivalently applied to the symmetrically split preconditioned (normal) system,
\begin{equation}
    (L^{-1})^\top \! B^\top \! BL^{-1}  (L\tilde{\sigma}) = (L^{-1})^\top \! B^\top  \tilde{b},
\end{equation}
where $L$ is (formally) a Cholesky factor of $H=L^\top \! L$. 

Let us finally return to the yet unspecified Dirichlet boundary condition that guarantees the positive definiteness of $H$. Since applying LSQR to \eqref{preconditioned} produces a sequence of approximate solutions that are in the range of $H^{-1}$ (see,~e.g.,~\cite{LSQR}) and  $\tilde{\sigma} = \sigma - \sigma_0$ is expected to vanish on average, it is natural to include in $H$ a homogeneous Dirichlet boundary condition on some subset of $\partial \Omega$. As in \cite{Edge-Enhancing}, we choose this subset to be the union of the electrodes $E$, which means that in the definition of $H$ in \eqref{H} the indices $i$ and $j$ only run over the basis functions that do not correspond to boundary nodes on the electrodes, whose number is defined to be $n - n' \in \mathbb{N}$. The degrees of freedom corresponding to these electrode nodes are always excluded from $\tilde{\sigma} \in \R^{n'}$ and included in $\sigma$ or $\sigma_0$ only when building the FE system matrix $A$ or for visualization purposes. Note that the rationale for implementing the boundary condition on $E$ is the hypothesis that the contact resistances can compensate for the effect of a possible misspecification of the conductivity level right underneath the electrodes \cite{Edge-Enhancing}. 

In order to consider the Morozov discrepancy principle as the stopping criterion for the LSQR iteration, observe that we are in essence solving the equation $B \tilde{\sigma} = \tilde{b}$ in a suitably regularized manner and, in addition,
\[
\tilde{b} - B \tilde{\sigma}  = b - B \sigma  =  \Gamma^{-1/2} P \big(y_P^{(i)} - J_{\sigma}^{(i)} \sigma  \big),
\]
where $P(y_P^{(i)} - J_{\sigma}^{(i)} \sigma)  = P \eta$ if the linearization and projection errors are ignored. Since $\Gamma^{-1/2}$ is the whitening matrix for $\eta$,
\[
 \mathbb{E} \big(\| \Gamma^{-1/2} P \eta \|^2 \big) \leq
\mathbb{E} \big( \| \Gamma^{-1/2} \eta \|^2 \big) = LM.
\]
In consequence, a reasonable choice for the noise level in the Morozov principle is 
\begin{equation}
\label{eq:morozov}
\epsilon = \sqrt{LM},
\end{equation}
i.e., the LSQR iteration is stopped when $\| b - B \tilde{\sigma} \| \leq \epsilon$. It is worth noting that if $\Gamma = \gamma^2 I$, as is the case in our numerical experiments, one can be more precise when computing the above expectation:
\[
 \mathbb{E} \big(\| \Gamma^{-1/2} P \eta \|^2 \big) =
\frac{1}{\gamma^2} \, \mathbb{E} \big( \| P \eta \|^2 \big) =  (L-1)M 
\]
because $P$ projects onto an $LM - M$ dimensional subspace and the variance of $\eta$ is the same in all directions. However, we anyway consider the mildly conservative noise level defined by \eqref{eq:morozov} in all our numerical experiments since the ratio $\sqrt{L/(L-1)}$ appears to be a reasonable fudge factor for avoiding a too low level of regularization.

The application of LSQR adds in principle a third nested iteration. However, our final reconstruction algorithm only takes one lagged diffusivity step for each sequential linearization, leaving only two nested loops, i.e., the sequential linearizations and LSQR, and making the index $j$ unnecessary. Based on a similar reasoning as above, a Morozov type stopping criterion is also used for the outer loop of sequential linearizations, 
\begin{equation}
     \big\lVert \Gamma^{-1/2}P(\mathbf{V}-\mathbf{U}(\sigma^{(i)}, z_0)) \big \rVert \leq \epsilon,
\end{equation}
where $\epsilon$ is as defined in \eqref{eq:morozov}.


\section{Implementation}\label{sec:implementation}

In this section, the numerical implementation is considered. We give enough details about the FE implementation to understand how the reduced basis method is applied and how it differs from the standard approach. In particular, we explain how $J_\sigma$ and $\mathbf{U}$ are computed at each step of the outer iteration and how this is accelerated with the help of a pre-computed reduced basis matrix $Q$. As the computation of $J_z = J_z(\sigma_0, z_0)$ happens during the offline stage, we do not consider its details but rather refer to, e.g., \cite{Vilhunen02}. To conclude the section, we present the projected lagged diffusivity iteration algorithm described in Section~\ref{sec:ProjectedLaggedDiffAlg} in the form in which it is numerically implemented.

\subsection{Computation of the Jacobians}

Recall from Section~\ref{sec:CEM} that the FE solution $(u,U)$ to \eqref{eq:cemeqs} is sought from the space $\{\psi_j\}_{j=1}^N \oplus {\rm span} \{\mathbf{c}^{(m)}\}_{m=1}^{M-1}$, where $\{\psi_j\}_{j=1}^N$ is a piecewise linear FE basis corresponding to a tetrahedral FE mesh of $\Omega$ and $\{\mathbf{c}^{(m)}\}_{m=1}^{M-1}$ is an orthonormal basis for $\R_\diamond^M$. We set
\begin{equation}
\mathcal{C}=\begin{bmatrix}
    \mathbf{c}^{(1)},\hdots , \mathbf{c}^{(M-1)}
\end{bmatrix} \in \R^{M \times (M-1)}
\end{equation}
and recall that the matrices $\mathcal{I} \in \R^{M \times L}$ and $\mathcal{U} (\sigma,z) \in \R^{M \times L}$ carry the employed current patterns and the corresponding computed electrode potentials, respectively, as their columns. The FE system matrix is still denoted by $A\in \R^{(N+M-1)\times (N+M-1)}$.

The nodal values of the FE solutions over all employed current patterns can be formally expressed as $A^{-1}f \in  \R^{(N+M-1) \times L}$ where $f=\begin{bmatrix}0 & (\mathcal{C}^\top \mathcal{I})^\top \end{bmatrix}^\top $. More precisely, the first $N$ elements in the $l$th column of $A^{-1}f$ give the nodal values of the interior potential $u$ for $I = I^{(l)}$, and the remaining  $M-1$ elements are the expansion coefficients of the associated electrode potentials in the basis $\{\mathbf{c}^{(m)}\}_{m=1}^{M-1}$. In particular,
\begin{equation}\label{Ubold}
   \mathcal{U}= \begin{bmatrix}0 & \mathcal{C} \end{bmatrix} A^{-1} \begin{bmatrix}0 \\[1mm] \mathcal{C}^\top \mathcal{I} \end{bmatrix}\in \R^{M\times L}.
\end{equation}
As quantified by \eqref{eq:U_to_U}, the vectorized forward solution $\mathbf{U}\in \R^{L M}$, extensively utilized in Section~\ref{sec:ProjectedLaggedDiffAlg}, is obtained from $\mathcal{U}$  by stacking its columns into a single vector. The not-yet-vectorized columns of the Jacobian $J_\sigma$ are thus obtained by differentiating \eqref{Ubold} with respect to $\sigma_l$,
\begin{equation}\label{Jk}
   (J_\sigma)_l = \frac{\partial \mathcal{U} }{\partial \sigma_l} = \begin{bmatrix}0 & \mathcal{C} \end{bmatrix} \frac{\partial A^{-1}}{\partial \sigma_l}  \begin{bmatrix}0 \\[1mm] \mathcal{C}^\top \mathcal{I} \end{bmatrix}= -\begin{bmatrix}0 & \mathcal{C} \end{bmatrix} A^{-1} \frac{\partial A}{\partial \sigma_l}  A^{-1}\begin{bmatrix}0 \\[1mm] \mathcal{C}^\top \mathcal{I} \end{bmatrix} \in \R^{M \times L},
\end{equation}
where $l=1, \dots, n$.

The stiffness matrix $A$ has an affine parameter dependence on $\sigma$. That is, it can be written in the form $A=A_0+\sum_{l=1}^n \sigma_l \delta A_l$, where $A_0$ and $\delta A_l$ do not depend on $\sigma$. Hence,
\begin{equation}
    \frac{\partial A}{\partial \sigma_l} = \delta A_l, \qquad l = 1, \dots, n,
\end{equation}
can be precomputed during the offline stage to speed up the formation of the Jacobian. The columns of $J_\sigma$ are then the vectorized versions of $(J_\sigma)_l$,
\begin{equation}\label{Jsigma}   
        J_\sigma = 
            \begin{bmatrix} 
        {\it vec}\big( (J_\sigma)_1 \big), \dots ,   
        {\it vec}\big( (J_\sigma)_n \big)
    \end{bmatrix} 
    \in \R^{LM \times n}.
\end{equation}

Although the number of electrodes $M$, and hence the size of $(J_\sigma)_l$, is typically small, evaluating \eqref{Jk} can still be costly due to the large size of the stiffness matrix~$A$ --- even more so considering that $n$ of these computations are needed at each sequential linearization in order to form the whole Jacobian in \eqref{Jsigma} at the considered iterate $\sigma^{(i)}$.


\subsection{Application of the reduced basis}

The use of $A^{-1}$ in the above is purely for notational simplicity. Explicitly inverting the large matrix $A$ would naturally be completely impractical, and one should instead resort to solving (sparse) systems of equations, which is in our case done with the backslash operator in MATLAB. Computing $\mathbf{U}$ and $J_\sigma$ thus amounts to first solving the auxiliary systems
\begin{equation}\label{AXAX0}
    AX = \begin{bmatrix}0 \\[1mm] \mathcal{C}^\top \mathcal{I} \end{bmatrix} \qquad \text{and} \qquad AX_0 = \begin{bmatrix}0 \\[1mm] \mathcal{C}^\top \end{bmatrix}
\end{equation}
for $X \in \R^{(N+M-1) \times L}$ and $X_0 \in \R^{(N+M-1) \times M}$, and then setting 
\begin{equation}
\label{eq:U_J_comput}
\mathbf{U} =  {\it vec}\big(\begin{bmatrix}0 & \mathcal{C} \end{bmatrix} \! X \big) \qquad \text{and} \qquad  (J_\sigma)_l = -X_0^\top \delta A_l X,
\end{equation}
where running through the indices $l=1, \dots, n$ determines the Jacobian $J_\sigma$.

This procedure can be accelerated by replacing $A$ with the much smaller reduced order surrogate $\hat{A}=Q^\top \!  AQ \in \R^{(k + M-1) \times (k + M-1)}$, where $Q \in \R^{(N+M-1) \times (k + M-1)}$, with $k \ll N$, is as described in Section~\ref{sec:RB}. Then instead of \eqref{AXAX0}, the reduced systems
\begin{equation}\label{hatAXAX0}
    \hat{A}\hat{X} = Q^\top \! \begin{bmatrix}0 \\[1mm] \mathcal{C}^\top \mathcal{I} \end{bmatrix} \qquad \text{and} \qquad \hat{A}\hat{X}_0 = Q^\top \! \begin{bmatrix}0 \\[1mm] \mathcal{C}^\top\end{bmatrix}
\end{equation}
are solved for $\hat{X} \in \R^{(k+M-1) \times L}$ and $\hat{X}_0 \in \R^{(k+M-1) \times M}$, and the redefined $X=Q\hat{X}$ and $X_0=Q\hat{X}_0$ are subsequently used in \eqref{eq:U_J_comput} when (approximately) evaluating $\mathbf{U}$ and $J_\sigma$. While $A$ is sparse, $\hat{A}$ is dense but has a much smaller dimension. In Section~\ref{sec:numerics}, we demonstrate numerically that such reduced order modeling leads to a significant speed-up with little effect on the accuracy of the reconstruction algorithm.

\begin{remark}
    The reduced basis version includes six additional multiplications by $Q^\top$ or $Q$. However, the extra cost of these operations is small compared to the time saved when solving the reduced systems \eqref{hatAXAX0} with $\hat{A}$ instead of \eqref{AXAX0} with the original system matrix~$A$.
\end{remark}

\begin{remark}
    The reduced basis is {\em not} employed when computing the preconditioner $H(\sigma)$ that includes information about the quick variations in $\sigma$. This is essential for being able to combine edge-enhancing regularization with reduced basis modeling. 
\end{remark}

\subsection{Algorithm}

The quantities that are updated at each outer iteration step, i.e., the quantitities that depend of the iterate $\sigma^{(i)}$, are marked with a raised parenthetical index $i$ to distinguish them from quantities that are available from the start or computed only once. There is no need for the index $j$ from Section~\ref{sec:LD} because we take only one lagged diffusivity step for each linearization. Since it is not the focus of this work, we refer to \cite{LSQR} for details on how to implement the LSQR iteration. We assume to have a good initial estimate for the background conductivity $\sigma_0$ as well as a homogeneous estimate for the contact resistances $z_0$. The reduced basis matrix $Q$ has been built during the offline stage, and its construction is not included in (the timing of) the algorithm presented below. Note that the initialization stage as well as building the projection operator can also be included in the offline phase.

The difference between the reduced basis and the standard implementation lies in whether the full or the reduced systems are solved in step 5. Note that this indirectly speeds up the computation of the Jacobian $J_\sigma$ in step 1  that uses the quantities computed in step~5.

\begin{algoritmi}\label{alg:ProjectedLaggedDiffusivity} Sequential linearization with projected lagged diffusivity step and LSQR
    \begin{itemize}
        \item[] \textbf{Initialization:} Choose $T>0$. Set ${\sigma}^{(0)}=\sigma_0$. Define the tolerances $\epsilon = \sqrt{LM}$ and $\delta > 0$.  Compute the FE system matrix at the initial guess $(\sigma_0, z_0)$. Solve for $X^{(0)}$ and $X_0^{(0)}$ (with or without reduced basis) and set $\mathbf{U}^{(0)}={ vec}(\begin{bmatrix}0 & \mathcal{C} \end{bmatrix} \! X^{(0)} )$. Set $i=0$.
        \item[] \textbf{Projection operator:} Compute $J_z$, the Jacobian matrix with respect to $z$ of the map $({\sigma}, z) \mapsto \mathcal{U}({\sigma}, z)$, evaluated at $(\sigma_0, z_0)$. Create $P = \mathrm{I} - J_z (J_z^\top J_z)^{-1} J_z^\top$, the orthogonal projection onto $\mathcal{R}(J_z)^\perp$. Define $\Sigma=\Gamma^{-1/2}P$.
    
    \item[] \textbf{Iteration:}
    \begin{enumerate}
        \item Using $X^{(i)}$ and $X_0^{(i)}$, compute $J_\sigma^{(i)}$,~i.e.,~the Jacobian with respect to ${\sigma}$ evaluated at $({\sigma}^{(i)}, z_0)$.
        \item Update $B^{(i)} = \Sigma J_\sigma^{(i)}$, $\tilde{\sigma}^{(i)}={\sigma}^{(i)}-\sigma_0$, $\tilde{y}_P^{(i)}=\mathbf{V}-\mathbf{U}^{(i)}+J_\sigma^{(i)} \tilde{\sigma}^{(i)}$, $\tilde{b}^{(i)}=\Sigma \tilde{y}_P^{(i)}$.
        \item Build the preconditioner $H(\tilde{\sigma}^{(i)})$ with a homogeneous Dirichlet boundary condition on $E$ for the purpose of operating with $H(\tilde{\sigma}^{(i)})^{-1}$.
        \item Apply the LSQR algorithm from \cite{LSQR} to (\ref{preconditioned}), stop according to the Morozov principle with the noise level $\epsilon$, and revert the change of variables to get ${\sigma}^{(i+1)}$. Replace negative components in ${\sigma}^{(i+1)}$ by $\delta$.
        \item Recompute the $\sigma$-dependent part of the FE system matrix for $\sigma^{(i+1)}$ with the nodal values on $E$ included.
        \item Forward solve for $X^{(i+1)}$ and $X_0^{(i+1)}$ (with or without reduced basis) and update $\mathbf{U}^{(i+1)}={\it vec}(\begin{bmatrix}0 & \mathcal{C} \end{bmatrix} \! X^{(i+1)}) $.
        \item Check the stopping criterion. Exit if $\lVert \Sigma (\mathbf{V}-\mathbf{U}^{(i+1)}) \rVert \leq \epsilon$. Otherwise, increment $i\rightarrow i+1$ and go to Step 1.
    \end{enumerate}
    \end{itemize}
\end{algoritmi}

\section{Numerical experiments}
\label{sec:numerics}

We test Algorithm~\ref{alg:ProjectedLaggedDiffusivity} in three cases that were considered in \cite{Edge-Enhancing} (Sections~\ref{sec:num1}--\ref{sec:num3}), as well as on data from another water tank (Section~\ref{sec:num4}). We refer to \cite{Edge-Enhancing} for the details of the setups of the first three experiments and mainly focus here on comparing reconstructions with and without reduced bases. Note that \cite{Edge-Enhancing} also considered the Perona--Malik prior that typically promotes edges more than the TV prior, which means that all reconstructions in \cite{Edge-Enhancing} are not directly comparable with the ones presented here.

The free parameters controlling the smoothing of the absolute value function in the definition of the TV prior and the lowest allowed conductivity level were, respectively, fixed to $T = 10^{-6}$ and $\delta = 10^{-2}$ in all experiments; according to our experience, moderate deviations from these values do not considerably affect the performance of Algorithm~\ref{alg:ProjectedLaggedDiffusivity}. In all experiments, the covariance of the additive zero-mean Gaussian noise process was assumed to be of the form
\begin{equation}
\label{eq:noise_cov}
\Gamma = \gamma^2 \mathrm{I} \in \R^{LM \times LM},
\end{equation}
where the common variance $\gamma^2$ was chosen in accordance with the specifications of the added noise for simulated data in Sections~\ref{sec:num1}--\ref{sec:num2} and based on the expected noise level for the considered measurement device --- as well as ``trial and error'' --- for experimental data in Sections~\ref{sec:num3}--\ref{sec:num4}.

When building reduced bases for the interior potential as described in Section~\ref{sec:RB}, the snapshot conductivites and contact resistances were drawn from log-normal distributions. More precisely, 
\begin{equation}
\label{eq:sample_sigma}
\log \sigma \sim \mathcal{N}(\log \sigma_0, \Gamma_0), 
\end{equation}
where the natural logarithm is applied componentwise, $\sigma_0 \in \R^n$ is the expected conductivity level of the target, and the covariance matrix is given elementwise as
\begin{equation} 
\label{eq:Gamma0}
(\Gamma_0)_{jk} = \omega^2 \exp \! \bigg ( \frac{\| x_j - x_k \|^2}{2 \ell^2} \bigg), \qquad j,k = 1, \dots, n.
\end{equation}
Here, $x_j$ denotes a node in the FE mesh employed for the discretization of the conductivity, and the pointwise standard deviation and the correlation length are, respectively, set to $\omega = 0.5$ and $\ell = 1$, which compare differently to the size of the imaged object in different experiments. Similarly, the contact resistances were assumed to be distributed according to
\begin{equation}
\label{eq:sample_z}
\log z \sim \mathcal{N}(\log z_0, \eta^2 \mathrm{I}), 
\end{equation}
where $z_0 = \zeta_0 \begin{bmatrix} 1, \dots, 1 \end{bmatrix}^\top$ represents the expected contact resistances, and the componentwise standard deviation was set to $\eta = 5 \cdot 10^{-4}$ in all numerical experiments. In particular, the expectations of the randomized versions of $\sigma$ and $z$, defined by \eqref{eq:sample_sigma} and \eqref{eq:sample_z}, do not exactly coincide with $\sigma_0$ and $z_0$, meaning that we intentionally introduced slight mismodeling when building the reduced basis. 

In all four experiments, the number of random draws was $K=500$, and the number of reduced basis vectors,~i.e.,~of the columns for the matrix $\widehat{Q}$ in \eqref{eq:tildeQ}, was $k=300$. The choice of $k$ was not optimized; it was simply deemed large enough for the reduced basis version of  Algorithm~\ref{alg:ProjectedLaggedDiffusivity} to produce reconstructions that are comparable in quality to those based on its standard version. On the other hand, $K$ was chosen large enough to ensure that the range of $\widehat{Q}$ was essentially not affected by the particular realizations of the random draws. Optimizing the size of $k$ would presumably lead to somewhat shorter run times for the reduced basis version of Algorithm~\ref{alg:ProjectedLaggedDiffusivity} without compromizing its reconstruction quality, whereas optimizing $K$ would allow a shorter offline phase for the algorithm.

As mentioned above, we do not account for the time required for building the reduced basis when comparing the standard and reduced basis implementations of Algorithm~\ref{alg:ProjectedLaggedDiffusivity}. The underlying assumption is that the reduced basis can be formed offline, that is, before having access to the measurement data. This requires knowing the measurement geometry as well as approximate conductivity and contact resistance levels for the target in advance. Although forming the reduced basis only requires a few tens of minutes of computation time for each of our numerical experiments, the reduced basis version of the algorithm is not competitive if this can only be done after having the measurements in hand. 

All computations were performed in a MATLAB environment on a laptop with 32 GB RAM and an Intel Core Ultra 5 125U CPU with maximum clock speed 4.3 GHz.

\subsection{Case 1: Cylinder (simulated data)}
\label{sec:num1}

\begin{figure}[t]
    \centering
    \begin{subfigure}{.49\textwidth}
	\includegraphics[width=1\textwidth]{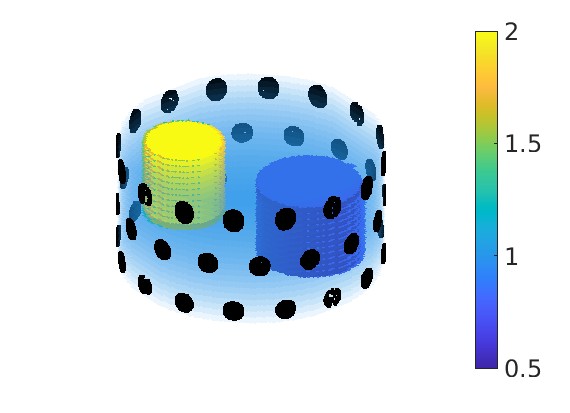}
    \caption{The target for Case~1.}
    \end{subfigure} 
    \begin{subfigure}{.49\textwidth}
    \includegraphics[width=1\linewidth]{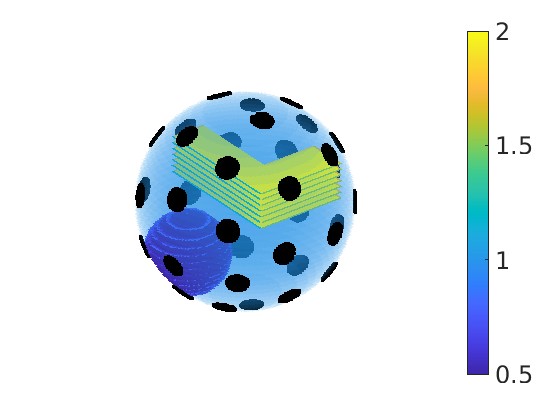}
    \caption{The target for Case~2.}
    \end{subfigure} 
    \caption{The measurement configurations and the target conductivities for Case~1 (left) and Case~2 (right).}
    \label{fig:cases}
\end{figure}

The setting of the first numerical experiment is visualized in the left image of Figure~\ref{fig:cases}. The examined object $\Omega$ is a cylinder of unit height and radius, with background conductivity level $\sigma = 1$ and $M = 48$ circular electrodes attached to its side. All electrodes are used for current input, meaning that the natural number of applied current patterns is $L=M-1$. Two cylindrical inclusions with conductivity levels $\sigma = 2$ and $\sigma = 0.5$ are embedded inside the object. Both have height $0.6$, with the former being slightly thinner and touching the top face of~$\Omega$, and the latter standing on the bottom face of $\Omega$. The contact resistances were drawn independently from a normal distribution with mean $2 \cdot 10^{-3}$ and standard deviation $5 \cdot 10^{-4}$, and the measurements were simulated on a dense FE mesh with $0.4\%$ of Gaussian noise added as a postprocessing step. For a more detailed description, see the first test case in \cite[Section~6]{Edge-Enhancing}.

\begin{figure}[t]
    \centering
        \begin{subfigure}{.45\textwidth}
	\includegraphics[width=1\textwidth]{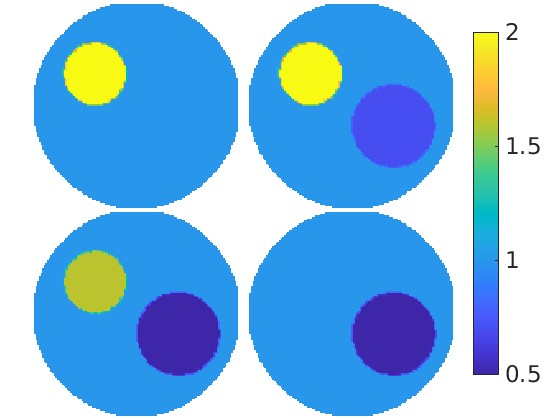}
    \caption{Target.}
    \end{subfigure} \quad
    \begin{subfigure}{.45\textwidth}
    \includegraphics[width=1\linewidth]{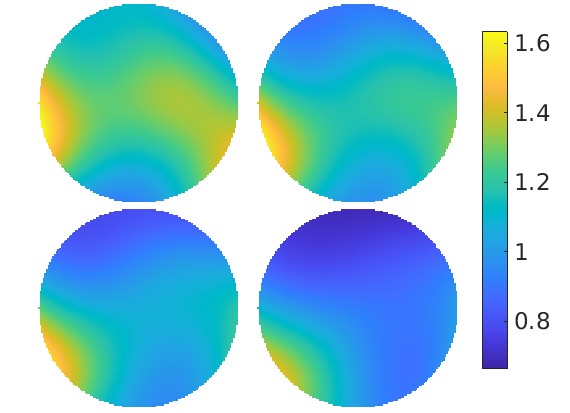}
    \caption{Example random draw from \eqref{eq:sample_sigma}.}
    \end{subfigure} \\[4mm]
    \begin{subfigure}{.45\textwidth}
		\includegraphics[width=1\textwidth]{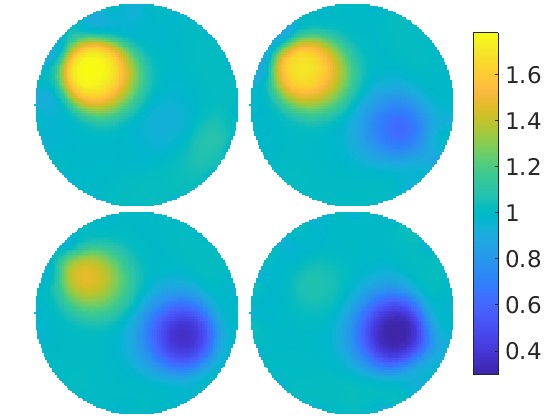}
        \caption{Reconstruction without reduced basis.}
    \end{subfigure} \quad
    \begin{subfigure}{.45\textwidth}
		\includegraphics[width=1\textwidth]{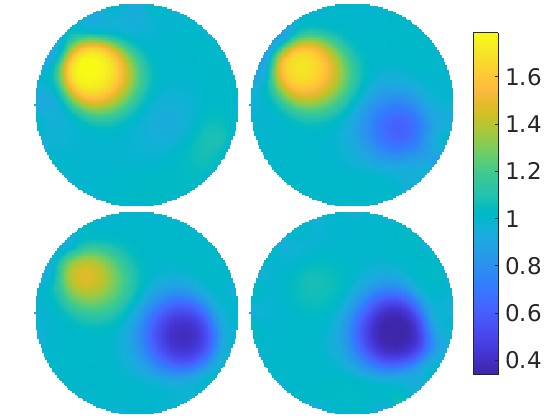}
        \caption{Reconstruction with reduced basis.}
    \end{subfigure}
    \caption{Case 1: The target is a cylinder of unit height and radius with two embedded cylindrical inclusions; see the left image of Figure~\ref{fig:cases}. The reconstruction and random draw cross sections are at heights $0.9$, $0.6$, $0.4$, and $0.1$.}
    \label{fig:case1}
\end{figure}

The FE mesh used for discretizing both the conductivity and the interior potential for Algorithm~\ref{alg:ProjectedLaggedDiffusivity} had $N = n =17\,031$ nodes and $80\,300$ tetrahedra. The employed estimates for the background conductivity and contact resistances were $\sigma_0 = 0.93$ and $\zeta_0 = 0.007$, respectively, and the other parameters in \eqref{eq:Gamma0} and \eqref{eq:sample_z} needed for forming the reduced basis for the interior potential were $\omega = 0.5$, $\ell = 1\,$, and $\eta = 5 \cdot 10^{-4}$. The top right image of Figure~\ref{fig:case1} presents an example draw from \eqref{eq:sample_sigma} --- notice that the correlation length for the log-normal density is in this case somewhat longer than the characteristic distances of the target inhomogeneities.

It took 43 seconds to run the standard version of Algorithm~\ref{alg:ProjectedLaggedDiffusivity} and 22 seconds to run its reduced basis version. As demonstrated by the bottom row of Figure~\ref{fig:case1}, which shows cross sections of the resulting two reconstructions, there is practically no difference in the reconstruction quality between the two algorithm versions. Moreover, the reconstructions are almost identical to the corresponding ones presented in \cite[Fig.~3]{Edge-Enhancing}; the run time reported in \cite{Edge-Enhancing} was 178 seconds with a CPU with clock speed 3.20 GHz, but these numbers are not directly comparable due to other differences in the hardware.

\subsection{Case 2: Ball (simulated data)}
\label{sec:num2}
In the second numerical experiment, $\Omega$ is a ball of radius 5 with constant background conductivity $\sigma = 1$, $M=37$ circular electrodes attached to its boundary, and two embedded inclusions. See the right image of Figure~\ref{fig:cases} for a visualization of the setup. One of the inclusions has conductivity $\sigma =0.5$ and is defined by the intersection of $\Omega$ with a ball of radius~2, whereas the other one is L-shaped with $\sigma = 2$. All electrodes are once again used for current input, leading naturally to $L = M-1$ linearly independent current patterns. The contact resistances were independently drawn from a normal distribution with mean  $10^{-2}$ and standard deviation $2 \cdot 10^{-3}$. The measurements were simulated on a dense FE mesh with $0.4\%$ of additive Gaussian noise. For further details, refer to the second test case in \cite[Section~6]{Edge-Enhancing}.

\begin{figure}[t]
    \centering
   \begin{subfigure}{.49\textwidth}
		\includegraphics[width=1\textwidth]{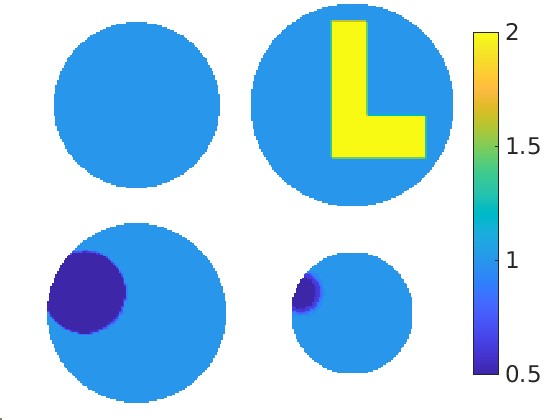}
        \caption{Target.}
    \end{subfigure} \ 
    \begin{subfigure}{.49\textwidth}
    \includegraphics[width=1\linewidth]{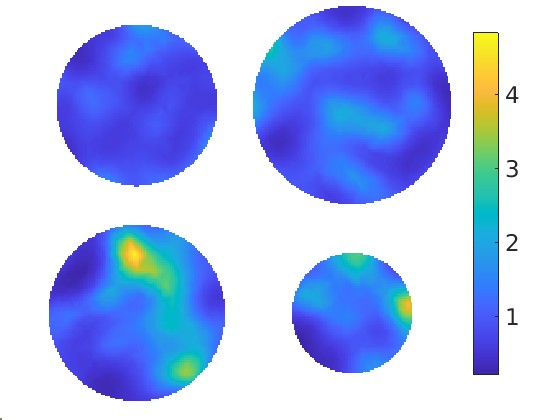}
    \caption{Example random draw from \eqref{eq:sample_sigma}.}
    \end{subfigure} \\[3mm]
    \begin{subfigure}{.49\textwidth}
		\includegraphics[width=1\textwidth]{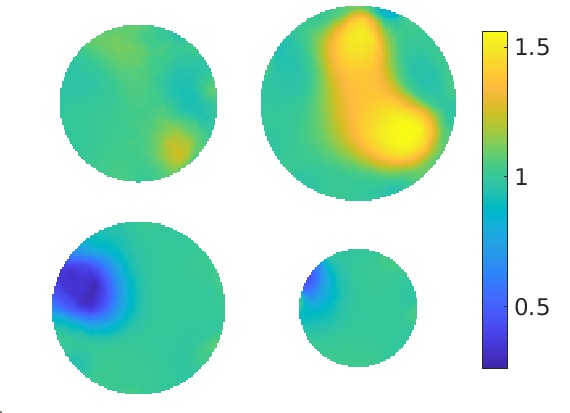}
        \caption{Reconstruction without reduced basis.}
    \end{subfigure}
    \ 
    \begin{subfigure}{.49\textwidth}
		\includegraphics[width=1\textwidth]{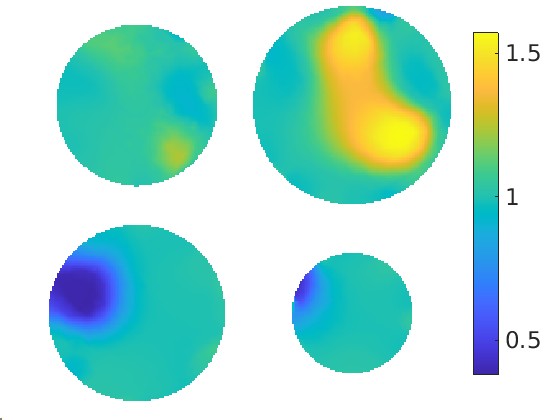}
        \caption{Reconstruction with reduced basis.}
    \end{subfigure}
    \caption{Case~2: The target is a ball of radius 5 with two embedded inclusions; see the right image of Figure~\ref{fig:cases}. The reconstruction and random draw cross sections are at heights $3$, $1$, $-2.5$, and $-4$.}
    \label{fig:case2}
\end{figure}

Algorithm~\ref{alg:ProjectedLaggedDiffusivity} employed a FE mesh with $N = n =11\,608$ nodes and $45\,896$ tetrahedra for both the interior potential and the conductivity reconstruction. The initial guesses for the background conductivity and contact resistances were, respectively, $\sigma_0 = 1.04$ and $\zeta_0 = 0.03$, and the other parameters in \eqref{eq:Gamma0} and \eqref{eq:sample_z} were chosen to be $\omega = 0.5$, $\ell = 1\,$, and $\eta = 5 \cdot 10^{-4}$. The top right image of Figure~\ref{fig:case2} depicts a random draw from \eqref{eq:sample_sigma} --- for this experiment, the correlation length for the log-normal density is well-aligned with the sizes of the target inclusions.

The runtime of the standard implementation of Algorithm~\ref{alg:ProjectedLaggedDiffusivity} was 15 seconds, which was cut down to 7 seconds for the reduced basis version. The bottom row of Figure~\ref{fig:case2} presents cross sections of the resulting two reconstructions that are almost exactly the same, as well as similar in quality with the corresponding one in \cite[Fig.~3]{Edge-Enhancing} for the Perona--Malik prior. The run time reported in \cite{Edge-Enhancing} was 86 seconds.

\subsection{Case 3: Thorax-shaped water tank (experimental data)}
\label{sec:num3}
\begin{figure}[t]
    \centering
    \begin{subfigure}{.40\textwidth}
		\includegraphics[width=1\textwidth]{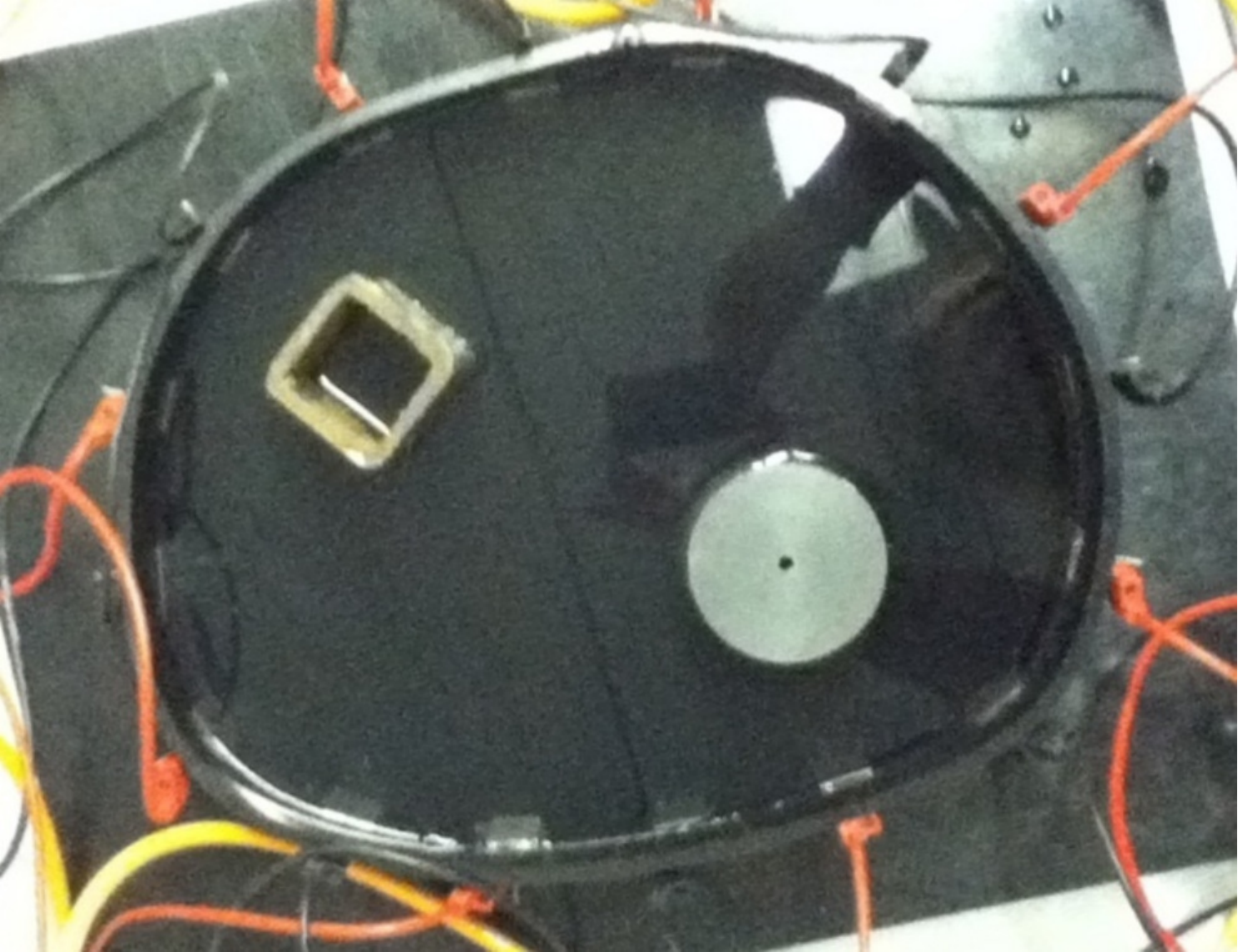}
        \caption{Measurement setup.}
    \end{subfigure} 
    \begin{subfigure}{.45\textwidth}
    \includegraphics[width=1\linewidth]{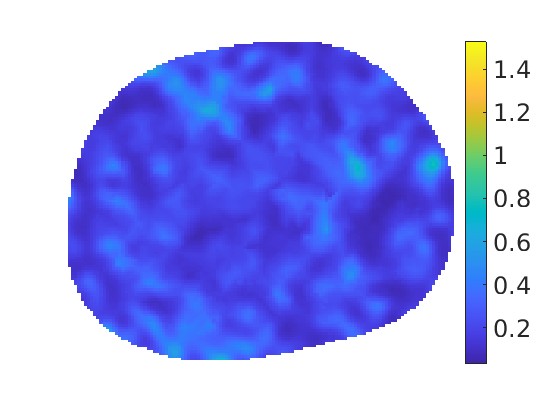}
    \caption{Example random draw from \eqref{eq:sample_sigma}.}
    \end{subfigure}
    \begin{subfigure}{.45\textwidth}
		\includegraphics[width=1\textwidth]{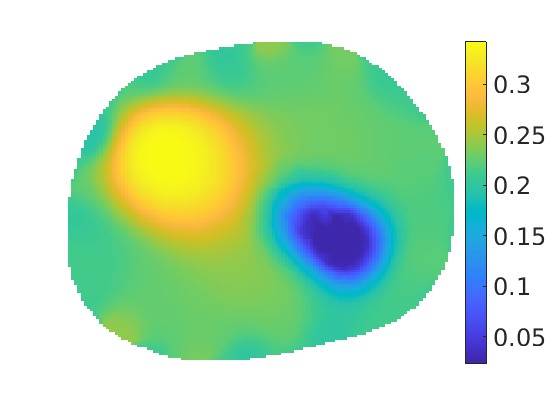}
        \caption{Reconstruction without reduced basis.}
    \end{subfigure}
    \begin{subfigure}{.45\textwidth}
		\includegraphics[width=1\textwidth]{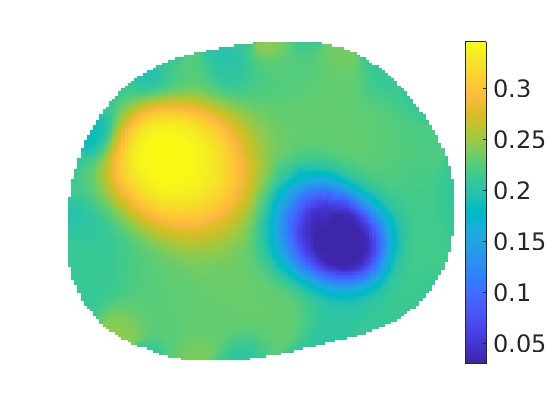}
        \caption{Reconstruction with reduced basis.}
    \end{subfigure}
    \caption{Case 3: Thorax-shaped water tank with two embedded cylindrical inclusions. The reconstruction and random draw cross sections are at height $2.5\,$cm and the unit of conductivity is mS/cm.}
    \label{fig:case3}
\end{figure}

The third test considers experimental data from a thorax-shaped water tank of circumference $106\,$cm filled with Finnish tap water and containing two inclusions: one made of insulating plastic and another made of highly conductive metal. The height of the water layer is $5\,$cm, and the inclusion cylinders extend from the bottom of the tank all the way through the water surface. There are $M=16$ rectangular electrodes attached to the side of the water tank, and $L = M-1$ linearly independent current patterns were fed in turns through the object. The measurements were performed with the KIT4 (Kuopio impedance tomography 4) device \cite{Kourunen09}. The measurement setup is shown in the top left image of Figure~\ref{fig:case3}; see the third test case in \cite[Section~6]{Edge-Enhancing} for further details.

The FE mesh used for the reconstruction algorithm had $N=n=27\,276$ nodes and $136\,286$ tetrahedra.
The noise covariance matrix was modeled according to \eqref{eq:noise_cov} with 
\begin{equation}
\label{eq:noise_std}
\gamma = \varsigma \max_{j,k} | \mathbf{V}_j - \mathbf{V}_k |,
\end{equation}
where we set $\varsigma = 6 \cdot 10^{-4}$ as in \cite{Edge-Enhancing}. The expected levels of conductivity and contact resistances were $\sigma_0 = 0.21\,$mS/cm and $\zeta_0 = 0.89\,{\rm k}\Omega \, {\rm cm}^2$. The other free parameters for forming the reduced basis were chosen as $\omega = 0.5$, $\ell = 1\,$cm, and $\eta = 5 \cdot 10^{-4}$, with the understanding that --- after applying the exponential function to the random draws --- the conductivity and contact resistances are given in the units mS/cm and ${\rm k}\Omega \, {\rm cm}^2$, respectively. A random draw from \eqref{eq:sample_sigma} is shown in the top right image of Figure~\ref{fig:case3}, according to which the employed correlation length for the log-normal density is much shorter than the diameters of the imaged inhomogeneities. 

Cross sections of the reconstructions produced by the standard and the reduced basis versions of Algorithm~\ref{alg:ProjectedLaggedDiffusivity} are shown on the bottom row of Figure~\ref{fig:case3}. Although there are some small differences between the reconstructions, especially close to the center of the tank, their quality and dynamic range are very similar. The run time was 17 seconds without and 7 seconds with the reduced basis. The corresponding reconstruction in \cite[Fig.~6]{Edge-Enhancing} has somewhat sharper edges, presumably due to the use of the Perona--Malik prior, but the overall reconstruction quality is arguably approximately the same. The runtime reported in \cite{Edge-Enhancing} was 180 seconds.

\subsection{Case 4: Circular water tank (experimental data)}
\label{sec:num4}

\begin{figure}[t]
  \centering
  \begin{subfigure}{.35\textwidth}
	\includegraphics[width=1\textwidth]{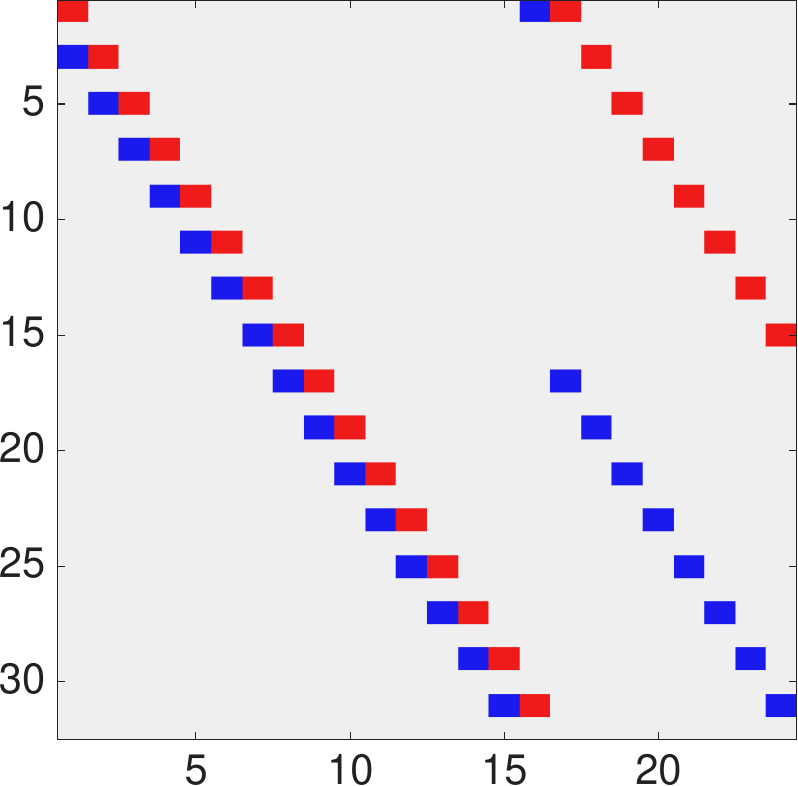}
        \caption{Current patterns.}
    \end{subfigure} 
    \quad 
    \begin{subfigure}{.5\textwidth}
\includegraphics[width=1\linewidth]{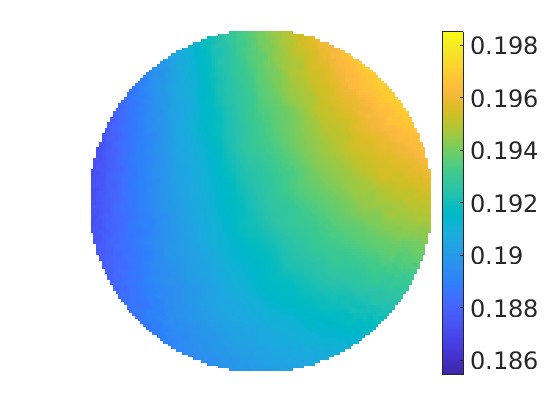}
    \caption{Example random draw from \eqref{eq:sample_sigma}.}
    \end{subfigure} 
   \caption{Case 4: Left: The $L = 24$ current patterns used for the cylindrical tank in Figure~\ref{fig:case4Afalse_and_true}. The vertical axis corresponds to electrode numbers and the horizontal axis to different current patterns. Red color indicates the source and blue color the sink of an 1\,mA current injection. Right: Random draw from \eqref{eq:sample_sigma}.}
    \label{fig:currents4}
\end{figure}

\begin{figure}[t]
   \centering
   \begin{subfigure}{.28\textwidth}
		\includegraphics[width=1\textwidth]{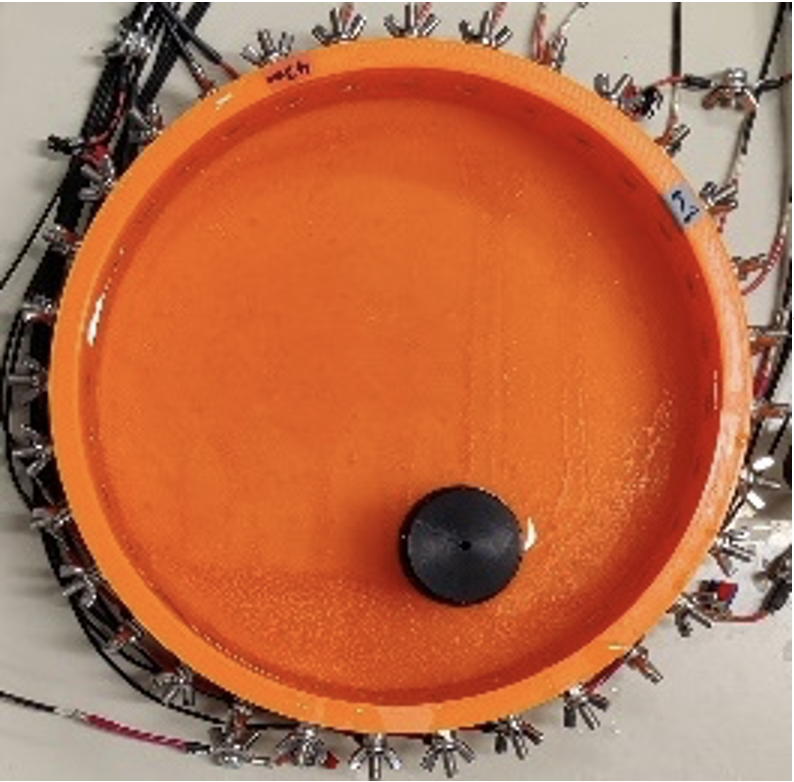}
        \caption{First target}
    \end{subfigure} \ \ \ \
    \begin{subfigure}{.29\textwidth}
\includegraphics[width=1\textwidth]{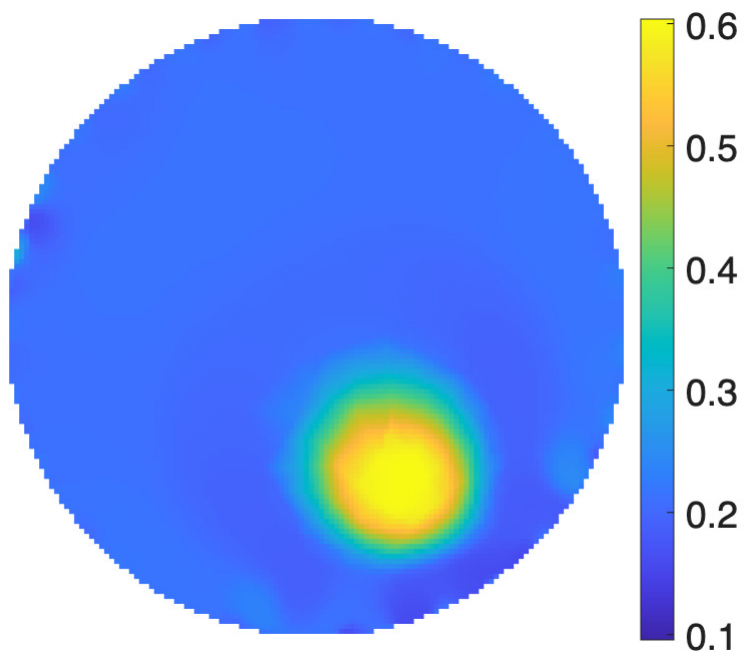}
        \caption{Reconstruction without reduced basis.}
    \end{subfigure} \ \ 
    \begin{subfigure}{.29\textwidth}
\includegraphics[width=1\textwidth]{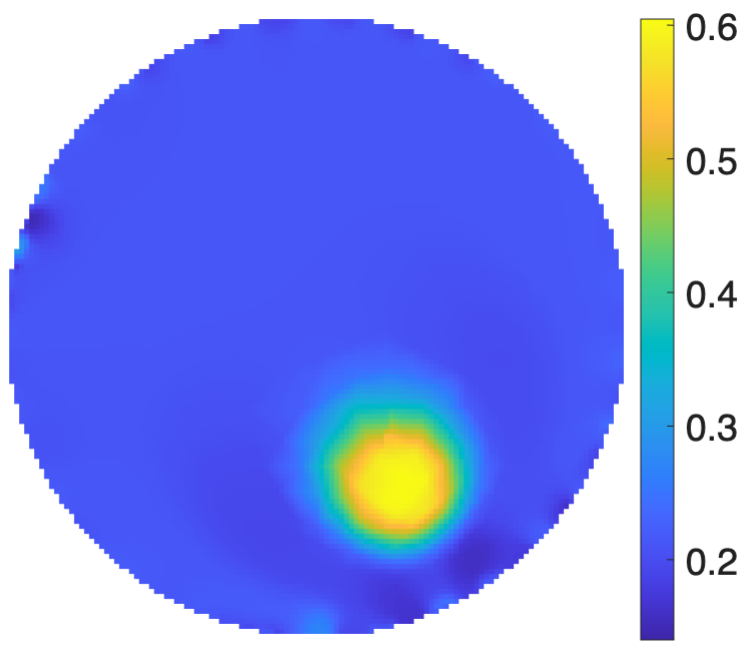}
        \caption{Reconstruction with reduced basis.}
    \end{subfigure}
    \\ [2mm]
    \begin{subfigure}{.28\textwidth}
		\includegraphics[width=1\textwidth]{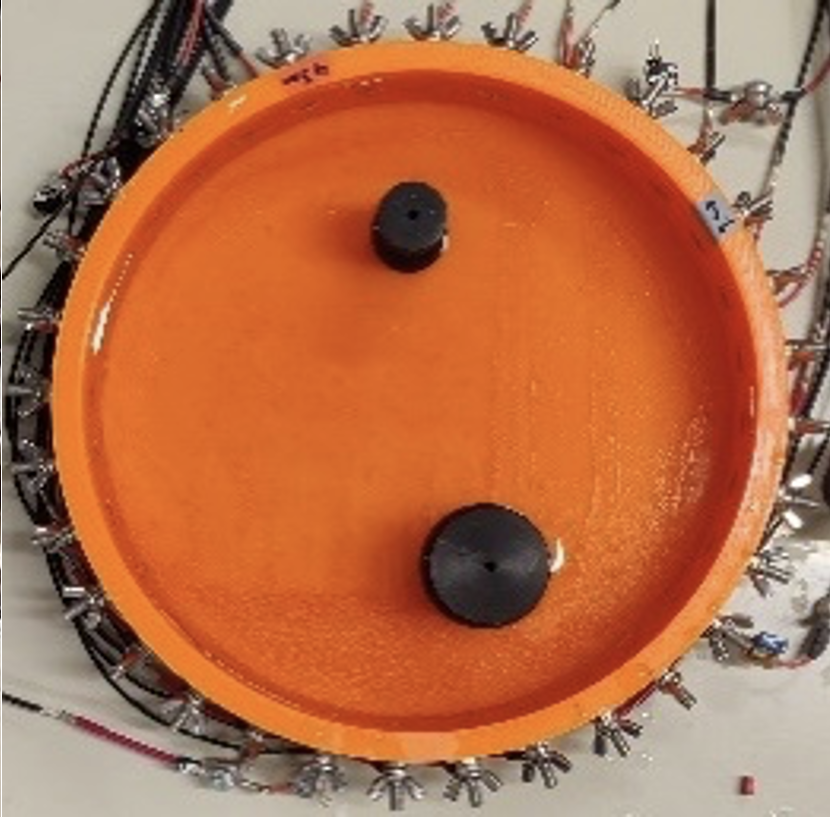}
        \caption{Second target.}
    \end{subfigure} \ \ \ \
    \begin{subfigure}{.29\textwidth}
\includegraphics[width=1\textwidth]{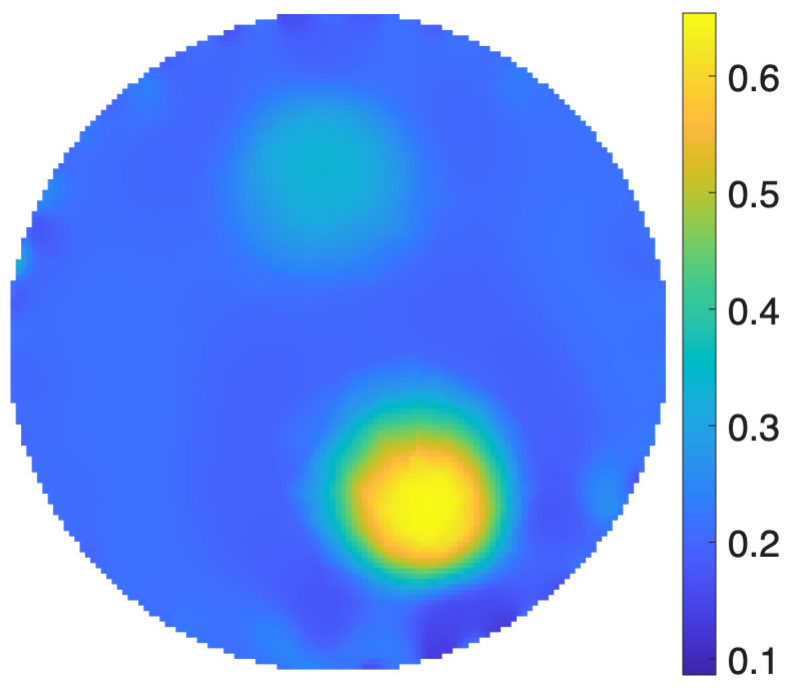}
        \caption{Reconstruction without reduced basis.}
    \end{subfigure} \ \ \
    \begin{subfigure}{.29\textwidth}
\includegraphics[width=1\textwidth]{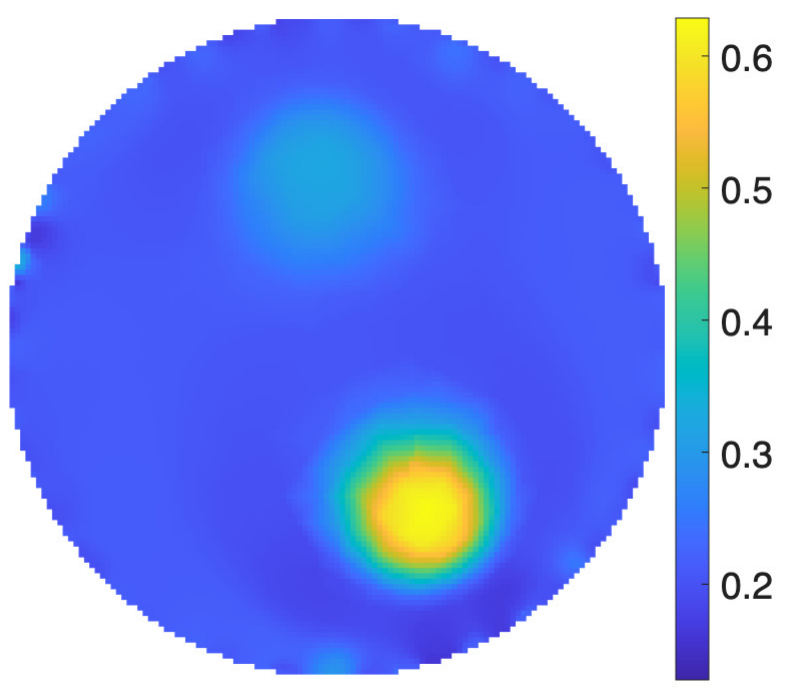}
        \caption{Reconstruction with reduced basis.}
    \end{subfigure}
    \caption{Case 4: Cylindrical water tank with either one (top) or two (bottom) conductive cylindrical inclusions. The displayed reconstruction cross sections are at height $2.15\,$cm, and the unit of conductivity is mS/cm.}
    \label{fig:case4Afalse_and_true}
\end{figure}

The fourth and final numerical experiment considers measurements on a cylindrical tank with radius 11.5\,cm shown on the left in Figure~\ref{fig:case4Afalse_and_true}. Depending on the test case, either one or two conductive cylindrical inclusion were placed inside the tank, and it was then filled with tap water to a depth of 4.3 cm. There are $M=32$ circular equiangularly positioned electrodes of radius 5\,mm attached to the tank's interior surface. The first electrode is identified by gray tape in Figure~\ref{fig:case4Afalse_and_true}, and the rest are numbered counterclockwise. The measurements were performed with the KIT5 (Kuopio impedance tomography 5) device described in detail in~\cite{Toivanen21}. The used $L=24$ current patterns are visualized in the left image of Figure~\ref{fig:currents4}; in particular, only the odd electrodes were used for driving currents. 

The FE mesh used for the reconstructions had $N=n=25\,412$ nodes and $122\,592$ tetrahedra. The noise covariance matrix was modeled according to \eqref{eq:noise_cov} and \eqref{eq:noise_std} with $\varsigma = 8 \cdot 10^{-5}$. The expected levels of conductivity and contact resistances were set to $\sigma_0 = 0.215\,$mS/cm and $\zeta_0 = 1\,{\rm k}\Omega \, {\rm cm}^2$; the former corresponds to the measured conductivity of the water in the tank, whereas the latter is simply a guess based on the contact resistances reconstructed in \cite{Edge-Enhancing}. In fact, by reconstructing a homogeneous estimate for the contact resistances as in \cite{Edge-Enhancing}, one would end up with a much lower value than $1\,{\rm k}\Omega \, {\rm cm}^2$. The reason for intentionally using a too high initial guess for the contact resistances is two-fold: it demonstrates that utilizing an orthogonal projection as explained in Section~\ref{sec:projection} indeed allows misspecifying the contacts without considerably affecting the reconstructions, and avoiding very low contact resistances in the computations reduces the effect that the singularities at the edges of the electrodes have on the numerical forward solutions,~cf.~\cite{Darde16}. The other free parameters for forming the reduced basis were chosen as $\omega = 0.5$, $\ell = 1$\,m, and $\eta = 5 \cdot 10^{-4}$. The right image in Figure~\ref{fig:currents4} visualizes a random draw from \eqref{eq:sample_sigma}, demonstrating that the reduced basis is built on a false assumption on the characteristic distances of the targets on the left in Figure~\ref{fig:case4Afalse_and_true}.

 The reconstructions produced by Algorithm~\ref{alg:ProjectedLaggedDiffusivity} without and with reduced basis are illustrated, respectively, in the center and right columns of Figure~\ref{fig:case4Afalse_and_true}. The two approaches once again result in reconstructions of similar quality, and the locations of the target inclusions are reproduced accurately in all four reconstructions. On the negative side, the reconstructions carry some small artifacts close to the boundary of the domain, possibly due to slight mismodeling of the measurement geometry. Moreover, the smaller cylindrical inclusion is not reconstructed accurately on the bottom row of Figure~\ref{fig:case4Afalse_and_true}: the reconstructed contrast is much lower than for the larger cylinder and the support of the smaller inclusion is also too extensive. The runtime for the standard version of Algorithm~\ref{alg:ProjectedLaggedDiffusivity} was 34.2 seconds for one inclusion and 35.5 seconds for two inclusions, while those for the reduced basis version were only 11.9 seconds and 12.6 seconds, respectively.

\section{Concluding remarks}
\label{sec:conclusion}
This work combined a reduced basis approach with a TV type prior in realistic EIT. More precisely, a POD was built for the forward solutions of the CEM in an offline phase,~i.e.,~prior to having access to the measurements but when already knowing the measurement geometry and some prior information on the expected conductivity and contact resistance levels for the target. Our numerical experiments on both simulated and experimental data demonstrated that the employment of reduced basis can accelerate the online phase of a TV-based reconstruction algorithm without negatively affecting the quality or edge-enhancing nature of the reconstructions. Although the reduced basis approach can speed up the computation of the forward solutions and Jacobian matrices required by an iterative reconstruction algorithm, the evaluation of the TV functional itself must still be performed with high spatial resolution, which restricts the achievable speed-up to a certain extent.

\section*{Acknowledgments}
We would like to thank Asko H\"anninen, Ville Kolehmainen and Jussi Toivanen from the University of Eastern Finland for the measurements used in the fourth numerical experiment. We would also like to express our gratitude to Altti J\"a\"skel\"ainen for his help with some details of the numerical implementation.

Vigdis Toresen's work was part of the Ministry of Education and Culture’s Doctoral Education Pilot under Decision No. VN/3137/2024-OKM-6 (Doctoral Education Pilot for Mathematics of Sensing, Imaging and Modelling). This work was also supported by the Finnish Research Council (decisions 353080, 353081, 358944, 359181).

\bibliographystyle{plain}
\bibliography{references}

\end{document}